\documentclass[a4paper,english,11pt,twoside,dvips]{amsart_schoemann}
\advance\oddsidemargin by -1.1cm
\advance\evensidemargin by -0.8cm
\usepackage{amsmath,amscd,amssymb,babel,epsfig,psfrag,haken}
\let\mathg\mathfrak
\theoremstyle{plain}
\newtheorem{thm}{Theorem}[section]
\newtheorem*{thm*}{Theorem}
\newtheorem{prop}{Proposition}[section]
\newtheorem{lem}{Lemma}[section]
\newtheorem{cor}{Corollary}[section]
\theoremstyle{definition}
\newtheorem{exa}{Example}[section]
\newtheorem{NB}{Remark}[section]
\newtheorem*{NB*}{Remark}
\newtheorem*{note*}{Note}

\newcommand{\bdm}{\begin{displaymath}}
\newcommand{\edm}{\end{displaymath}}
\newcommand{\x}{\times}
\newcommand{\op}{\oplus}
\newcommand{\ox}{\otimes}

\newcommand{\ra}{\rightarrow}

\newcommand{\lan}{\left\langle}
\newcommand{\ran}{\right\rangle}

\newcommand{\C}{\ensuremath{\mathbb{C}}}
\newcommand{\R}{\ensuremath{\mathbb{R}}}
\newcommand{\N}{\ensuremath{\mathbb{N}}}
\newcommand{\Z}{\ensuremath{\mathbb{Z}}}

\renewcommand{\u}{\ensuremath{\mathfrak{u}}}
\newcommand{\su}{\ensuremath{\mathg{su}}}
\newcommand{\so}{\ensuremath{\mathg{so}}}

\newcommand{\g}{\ensuremath{\mathg{g}}}
\newcommand{\h}{\ensuremath{\mathg{h}}}
\newcommand{\m}{\ensuremath{\mathg{m}}}
\newcommand{\n}{\ensuremath{\mathg{n}}}
\newcommand{\slla}{\ensuremath{\mathg{sl}}}
\newcommand{\hol}{\ensuremath{\mathg{hol}}}

\newcommand{\isomorph}{\cong}			% isomorph (rather than congruent)
\newcommand{\iso}{\mathg{iso}} 			% Lie algebra of isotropy group
\newcommand{\lsemidirect}{\ltimes}		% semi-direct product (left)
\newcommand{\hodge}{*}		                % hodge
\newcommand{\cycsum}{\mathfrak{S}}              % cyclic sum 
\newcommand{\Lie}[1]{\mcL_{#1}}                 % Lie derivative
\newcommand{\CP}{\mathbb{CP}}                   % CP(n)
\newcommand{\mcF}{\mathcal{F}}            
\newcommand{\mcH}{\mathcal{H}}            
            
\newcommand{\mcL}{\mathcal{L}}            
\newcommand{\mcR}{\mathcal{R}}            
            
\newcommand{\mcV}{\mathcal{V}}            

\newcommand{\spacedhline}{\\[1pt]\hline}

\let\LEFTMARK=\leftmark\let\RIGHTMARK=\rightmark
\def\rightmark{\smash{\RIGHTMARK}\vadjust{\vskip1.5mm\hrule height 0.4pt}}
\def\leftmark{\smash{\LEFTMARK}\vadjust{\vskip1.5mm\hrule height 0.4pt}}
\begin{document}
%%%%%%%%%%%%%%%%%%%%%%%%%%%%%%%%%%%%%%%%%%%%%%%%%%%%%%%%%%%%%%%%%%%%%%%%%%%%%
%
\title[Almost hermitian structures with parallel torsion]{Almost hermitian structures\\ with parallel torsion}
\author{Nils Schoemann}
\address[Schoemann]{Institut f\"ur Mathematik\\ 
Humboldt-Universit\"at zu Berlin\\
Unter den Linden 6\\
Sitz: John-von-Neumann Haus, Adlershof\\
D-10099 Berlin, Germany}
\email{schoeman@mathematik.hu-berlin.de}
\date{7th March 2007}
%\date{\today}
%
\subjclass[2000]{Primary 53 C 25; Secondary 81 T 30}
\thanks{Supported by the International Research Training Group "Arithmetic and Geometry" (GRK 870) of the German Research Foundation (DFG)}
\begin{abstract}
The characteristic connection of an almost hermitian structure is a hermitian connection with totally skew-symmetric torsion. The case of parallel torsion in dimension six is of particular interest. In this work, we give a full classification of the algebraic types of the torsion form, and, based on this, undertake a systematic investigation into the possible geometries. Numerous naturally reductive spaces are constructed and classified as well as examples on nilmanifolds given.
\end{abstract}
\maketitle
%
%
%\tableofcontents
%
\pagestyle{headings}
%
%\newpage
%
%%%%%%%%%%%%%%%%%%%%%%%%%%%%%%%%%%%%%%%%%%%%%%%%%%%%%%%%%%%%%%%%%%%%%%%%%%%%%
\section{Introduction}
%%%%%%%%%%%%%%%%%%%%%%%%%%%%%%%%%%%%%%%%%%%%%%%%%%%%%%%%%%%%%%%%%%%%%%%%%%%%%
\noindent
Making systematic use of the covariant derivative of the K\"ahler form, A. Gray and L. M. Hervella, in the late seventies, divided almost hermitian structures into sixteen classes. The class of cocalibrated structures, class $W_1\op W_3\op W_4$, is interesting, for precisely those structures admit a hermitian connection with totally skew-symmetric torsion. The connection in that case is unique, and is called the characteristic connection of the almost hermitian structure. Examples of cocalibrated structures include nearly K\"ahler structures, class $W_1$, and hermitian structures, class $W_3\op W_4$. 

About the same time, V. F. Kirichenko proved, that the torsion form of the characteristic connection of any nearly K\"ahler structure is parallel. For naturally reductive spaces equipped with an invariant almost hermitian structure, the canonical and characteristic connection coincide. The general theory of homogenous spaces implies, that the torsion form is parallel. Further examples of structures with parallel torsion include generalized Hopf structures, introduced and first studied by I. Vaisman, and certain isolated, but important, examples on twistor spaces.

In this paper, we treat almost hermitian structures in dimension six such that the characteristic connection exists and the torsion form is parallel. Principally, we do the following:

\begin{enumerate}
\item[(i)] We study the geometric invariants associated with such structures: the isotropy group of the torsion form, the holonomy group of the characteristic connection, and the (precise) algebraic type of the torsion form.
\item[(ii)] We obtain an overview of all such structures, and study the geometry of certain classes in more detail.
\end{enumerate}

\noindent
After an introductory part, where we collect basic material on almost hermitian structures, we give a full classification of the algebraic type of the torsion form. Then, based on this, we undertake a systematic investigation into the possible geometries.

\begin{note*}
This work is a continuation of the work of B. Alexandrov, Th. Friedrich, and the author in \cite{AlexandrovFriedrichSchoemann:2005}.
\end{note*}

%\newpage
%%%%%%%%%%%%%%%%%%%%%%%%%%%%%%%%%%%%%%%%%%%%%%%%%%%%%%%%%%%%%%%%%%%%%%%%%%%%%
\section{Elements of almost hermitian geometry}
%%%%%%%%%%%%%%%%%%%%%%%%%%%%%%%%%%%%%%%%%%%%%%%%%%%%%%%%%%%%%%%%%%%%%%%%%%%%%
%
%---------------------------------------------------------------------------- 
\subsection{Definitions} 
%---------------------------------------------------------------------------- 
Let $(M^6,g,J)$ denote an almost hermitian structure. That is, a smooth manifold $M^6$ together with (positive definite) Riemannian metric $g$, and an orthogonal almost complex structure $J$. The \emph{K\"ahler form} is the real valued skew-symmetric form defined by
\begin{gather*}
\Omega(X,Y)=g(JX,Y)
\end{gather*}
for any vectors $X,Y$. A frame $e=(e_1,...,e_6)$ is positively oriented (or simply oriented) if $(\Omega\wedge\Omega\wedge\Omega)(e_1,...,e_6)>0$. There exists an oriented orthonormal frame such that
\begin{gather*}
\Omega=e_{12}+e_{34}+e_{56},
\end{gather*}
and any such frame is called an \emph{adapted frame}. (Here and henceforth we do not distinguish notationally between vectors and covectors, and we use the short hand $e_{ij}$ for the exterior product $e_i\wedge e_j$, and so on.) Let $\mcF$ and $\mcR$ denote the bundle of oriented orthonormal and adapted frames, with structure groups $SO(6)$ and $U(3)$ respectively.

%---------------------------------------------------------------------------- 
\subsection{Modules}
%---------------------------------------------------------------------------- 
Let $\R^6=(\R^6,g,J)$ be the standard hermitian vector space, $\Omega$ the K\"ahler form, and $(e_1,...,e_6)$ an adapted basis.

The module of 2-forms $\Lambda^2(\R^6)$ is isomorphic to the module of skew-symmetric endomorphisms $\so(6)=\so(\R^6)$. Let $\u(3)$ and $\m^6$ be the submodules of endomorphisms that commute respectively anti-commute with the complex structure, 
\begin{gather*}
\u(3)=\{A\in\so(6):AJ=JA\},\quad\m^6=\{A\in\so(6):AJ=-JA\}.
\end{gather*}
Let $\su(3)\subset\u(3)$ be the trace free part (complex trace), and $\R^1\subset\u(3)$ the submodule spanned by the complex structure.

\begin{prop}
As a $U(3)$-module,
\begin{gather*}
\so(6)=\u(3)\op\m^6=\R^1\op\su(3)\op\m^6.
\end{gather*}
The $U(3)$-modules $\su(3)$ and $\m^6$ are real irreducible. The module $\m^6$ is isomorphic to $\C^3$ with action $a\cdot v=det(a)^3av$ for $a\in U(3)$, $v\in\C^3$. 
\end{prop}

The module of 3-forms $\Lambda^3(\R^6)$ is a complex module (the restriction of the Hodge operator is a complex structure). The module is complex reducible. To describe the irreducible components consider the operator $\tau:\Lambda^3(\R^6)\ra\Lambda^3(\R^6)$ obtained by evaluation of the differential of the representation $\rho:U(3)\ra Aut(\Lambda^3(\R^6))$ on the complex structure, $\tau=\rho_*(J)$. Explicitly, $\tau(T) = \sum_{i=1}^6(e_i\haken\Omega)\wedge(e_i\haken T)$. The operator $\tau$ is $U(3)$-equivariant and commutes with the complex structure. Its square $\tau^2$ has eigenvalues -9, -1 and 1. The corresponding eigenspaces are of real dimension 2, 12 and 6, and are denoted by $\Lambda^3_2$, $\Lambda^3_{12}$ and $\Lambda^3_6$. 

\begin{prop}\label{prop_module_Lambda3}
As a $U(3)$-module,
\begin{gather*}
\Lambda^3(\R^6)=\Lambda^3_2\op\Lambda^3_{12}\op\Lambda^3_6.
\end{gather*}
The $U(3)$-modules $\Lambda^3_2$, $\Lambda^3_{12}$ and $\Lambda^3_6$ are real irreducible. The module $\Lambda^3_2$ is isomorphic to $\C$ with action $a\cdot z=det(a)^3z$ for $a\in U(3)$, $z\in\C$, and, in particular, is a trivial $SU(3)$-module. The module $\Lambda^3_6$ is isomorphic to $\R^6=\C^3$, the standard $U(3)$-module. The map $\R^6\ra\Lambda^3_6$ with $X\mapsto\Omega\wedge X$ is an isomorphism. 
\end{prop}

\label{txt_module_R6xm6}
The module $\R^6\ox\m^6$ is the tensor product of two complex modules. It is complex reducible, and the irreducible submodules can be characterized in terms of the irreducible submodules of $\Lambda^3(\R^6)$. Let $\theta:\Lambda^3(\R^6)\ra\R^6\ox\m^6$ be the $U(3)$-equivariant map defined by
$\theta(T)=-\frac{1}{2}\sum_{i=1}^6e_i\ox pr_{\m^6}(e_i\haken T)$, where $pr_{\m^6}$ denotes the projection $\so(6)\ra\m^6$. It is injective. Let $W_1=\theta(\Lambda^3_2)$, $W_3=\theta(\Lambda^3_{12})$, $W_4=\theta(\Lambda^3_6)$, and let $W_2$ denote the orthogonal complement of $W_1\op W_3\op W_4$ inside $\R^6\ox\m^6$.
\begin{prop}
As a $U(3)$-module,
\begin{gather*}
\R^6\ox\m^6=W_1\op W_2\op W_3\op W_4.
\end{gather*}
The $U(3)$-modules $W_1$, $W_2$, $W_3$ and $W_4$ are real irreducible.
\end{prop}

%---------------------------------------------------------------------------- 
\subsection{Intrinsic torsion, the sixteen classes, and the characteristic connection}
%---------------------------------------------------------------------------- 
The intrinsic torsion of an almost hermitian structure measures the failure of the Levi-Civita connection to reduce from a connection in the bundle of oriented orthonormal frames to connection in the bundle of adapted frames. The connection reduces if and only if the structure is K\"ahler.

Let $\omega$ be the connection 1-form of the Levi-Civita connection $\nabla^g$ seen as a function on $\mcF$ with values in $\R^6\ox\so(6)$, that is $\omega_{ij}(e)=g(\nabla^ge_i,e_j)$ for $e\in\mcF$. Let $\omega$ also denote the restriction from $\mcF$ to $\mcR$. Let $pr_{\R^6\ox\m^6}$ be the projection $\R^6\ox\so(6)\ra\R^6\ox\m^6$. The section of the bundle $\mcR\x_{U(3)}\R^6\ox\m^6$ defined by $\Gamma=pr_{\R^6\ox\m^6}(\omega)$ is called the \emph{intrinsic torsion} of the almost hermitian structure. The decomposition of the $U(3)$-module $\R^6\ox\m^6$ into irreducible submodules determines four sections, $\Gamma=\Gamma_2+\Gamma_{16}+\Gamma_{12}+\Gamma_6$. Here, $\Gamma_i$ is the projection of $\Gamma$ onto one of the submodules of $\R^6\ox\m^6$, the subscript denoting the real dimension of that module. 

The (pointwise) invariants type and strict type are defined corresponding to the vanishing or non-vanishing of the various components of the intrinsic torsion. The almost hermitian structure is said to be of \emph{type} $W$ if the intrinsic torsion takes values in $W$. It is said to be of \emph{strict type} $W$ if it is of type $W$ and not of type $W'$ for all submodules $W'\varsubsetneq W$. The four submodules $W_1$, $W_2$, $W_3$ and $W_4$ determine combinatorially a maximum of sixteen submodules. These are the \emph{sixteen classes} of almost hermitian structures. 

Let $T_2$, $T_{12}$, $T_6$ be the 3-forms determined by the components $\Gamma_2$, $\Gamma_{12}$ and $\Gamma_6$ of the intrinsic torsion, where the identification is made such that $\theta(T_i)=\Gamma_i$. Let $X$ be the vector field determined by $T_6=\Omega\wedge X$. (The map $\Omega\wedge..:\R^6\ra\Lambda^3_6$ is an isomorphism.) The 3-form and the vector field both will be used to describe the $W_4$-component of the intrinsic torsion.

A \emph{hermitian connection} is a connection $\nabla$ such that metric and almost complex structure are parallel, $\nabla g=0$, $\nabla J=0$. The torsion, defined as $(2,1)$-tensor by the formula $T(X,Y)=\nabla_XY-\nabla_YX-[X,Y]$ for vector fields $X,Y$, and seen as $(3,0)$-tensor through the formula $T(X,Y,Z)=g(T(X,Y),Z)$, is skew-symmet\-ric in $X$ and $Y$ by definition and is called  \emph{totally skew-symmetric} if it is skew-symmetric in $Y$ and $Z$. In that case, $T$ is called the \emph{torsion form}. 

The following theorem characterizes the various classes through sets of differential equations. 

\begin{thm}\label{thm_type_and_diff_eqs}
Let $(M^6,g,J)$ be an almost hermitian structure. Then differential and codifferential of the K\"ahler form satisfy
\begin{gather*}
d\Omega=3\hodge T_2-\hodge T_{12}+\hodge T_6,\quad\delta\Omega = 2X.
\end{gather*}
The Nijenhuis tensor is totally skew-symmetric if and only if the structure is of type $W_1\op W_3\op W_4$, and vanishes if and only if the structure is of type $W_3\op W_4$. There exists a hermitian connection with totally skew-symmetric torsion if and only if the structure is of type $W_1\op W_3\op W_4$. In that case the connection is unique, 
\begin{gather*}
\nabla^c=\nabla^g+\frac{1}{2}\;T^c,
\end{gather*} 
and its torsion form is $T^c=T_2+T_{12}+T_6$.
\end{thm}

\noindent
The connection is called the characteristic connection of the almost hermitian structures. The \emph{type} and \emph{strict type} of the torsion form are defined in obvious ways.\footnote{The Nijenhuis tensor of an almost complex structure is the (2,1)-tensor defined by $N(X,Y)=[JX,JY]-J[JX,Y]-J[X,JY]-[X,Y]$ for any vector fields $X,Y$. It is seen as a (3,0)-tensor through the formula $N(X,Y,Z)=g(N(X,Y),Z)$.}

%---------------------------------------------------------------------------- 
\subsection{Differentials and Curvature}
%---------------------------------------------------------------------------- 
From now on, we will consider almost hermitian structures such that the characteristic connection exists. 

Let $\alpha$ be a $k$-form which is parallel, $\nabla^c\alpha=0$. The differential then is given by
\begin{gather*}
d\alpha=\sum_{i=1}^{n}(e_i\haken\alpha)\wedge(e_i\haken T^c).
\end{gather*}
For 1-forms the formula $d\alpha=\alpha\haken T^c$ holds. For the torsion form itself, if parallel, we have $dT^c=2\sigma_{T^c}$, where here and throughout $\sigma_{T^c}=\frac{1}{2}\sum_{i=1}^{n}(e_i\haken T^c)\wedge(e_i\haken T^c)$.\label{page_sigmaTc}

The curvature of the characteristic and the Levi Civita connection are related through the formula
\begin{gather*}
(R^c-R^g)(X,Y,Z,U)=\;\frac{1}{2}(\nabla^c_XT^c)(Y,Z,U)-\frac{1}{2}(\nabla^c_YT^c)(X,Z,U) \\
+\frac{1}{4}g(T^c(X,Y),T^c(Z,U))+\frac{1}{4}\sigma_{T^c}(X,Y,Z,U),
\end{gather*}
for vector fields $X,Y,Z,U$. If the torsion is parallel, this difference is an algebraic expression in the torsion form (quadratic). Similar formulas hold for the Ricci and scalar curvature, as they are obtained by contraction. The first Bianchi identities\label{page_Biachi_identities} take the form
\begin{gather*}
(\cycsum R^c)(X,Y,Z,U)=\cycsum(\nabla_XT^c)(Y,Z,U)+\sigma_{T^c}(X,Y,Z,U).
\end{gather*}
Here $\cycsum$ denotes the cyclic sum in $X,Y,Z$. If the torsion is parallel, this is an an algebraic expression in the curvature tensor and the torsion form. 

\begin{NB}\label{rmk_Sternberg_criteria}
In the case of Lie groups, i.e. parallel torsion and trivial holonomy, the curvature vanishes and the Bianchi identities are a quadratic expression in the torsion form. The Bianchi identities hold if and only if the square of the torsion form inside the Clifford algebra is a scalar, $(T^c)^2\in\R\subset Cl(\R^6)$, (cf. \cite{Sternberg:2004}).
\end{NB}

%\newpage
%%%%%%%%%%%%%%%%%%%%%%%%%%%%%%%%%%%%%%%%%%%%%%%%%%%%%%%%%%%%%%%%%%%%%%%%%%%%%
\section{Torsion Orbits}
%%%%%%%%%%%%%%%%%%%%%%%%%%%%%%%%%%%%%%%%%%%%%%%%%%%%%%%%%%%%%%%%%%%%%%%%%%%%%
%
\noindent
The torsion form is a map $T^c:\mcR\ra\Lambda^3(\R^6)$, such that $T^c(ea)=a^{-1}T^c(e)$ for any $e\in\mcR$, $a\in U(3)$. At each point of the manifold, it determines an $U(3)$-orbit in $\Lambda^3(\R^6)$. If the torsion form is parallel the orbit is constant (on each connected component). It is natural to ask for the set of orbits that may occur as the orbit of the torsion form, the \emph{torsion orbits}. 

The \emph{isotropy group} of the orbit containing $T^c\in\Lambda^3(\R^6)$ is the $U(3)$-conjugacy class of the group $Iso(T^c)=\{a\in U(3):aT^c=T^c\}$. The \emph{restricted isotropy group} $Iso_o(T^c)$ is the connected component of the identity of the isotropy group. The geometric significance of the isotropy group derives from the fact, that, if the torsion form is parallel, the holonomy group is contained in the isotropy isotropy group,
$Hol(\nabla^c)\subset Iso(T^c)$. 

%----------------------------------------------------------------------------
\subsection{Main classification theorem}
%----------------------------------------------------------------------------
The conjugacy class of a subgroup $G\subset U(3)$, may be described by specifying an abstract group $G$, a faithful representation $G\ra Aut(\C^3)$, and a $G$-invariant hermitian form. The hermitian form here need not be unique, the choice of various hermitian forms may determine the same conjugacy class. In the following, if $G$ may be embedded in $U(3)$ explicitly, a representation and (some) hermitian form are understood.

The group $SU(3)$ is a maximal subgroup of $U(3)$, and all groups isomorphic to $SU(3)$ are conjugate. For the group $SO(3)$, note that the complexification of its standard representation on $\R^3$ yields its unique complex 3-dimensional representation. Since there is a unique invariant hermitian form, it follows that $U(3)$ contains precisely one conjugacy class of groups isomorphic to $SO(3)$. For $k\in\Z$, the embeddings of the group $U(2)$ defined by $U(2)_k=\{\bigl[\begin{smallmatrix} A & 0 \\ 0 & det(A)^k \end{smallmatrix}\bigl]: A\in U(2)\}$ determine pairwise non-conjugate subgroups. The groups $U(2)_0$ and $U(2)_1$ occur as the isotropy group of some torsion form. The group $U(2)_{-1}\subset SU(3)$ does not, but it does occur as holonomy group. For $k\neq-1,0,1$, the groups $U(2)_k$ do not preserve any non-trivial 3-form. Henceforth, let $SU(2)$ be the maximal subgroup of $U(2)_0$. Moreover, let $T^2$ be the maximal torus in $U(2)_0$, and $T^1$ the maximal torus in $SU(2)$.

\begin{note*}
There are many non-equivalent embeddings of abstract tori of dimension one and two into $U(3)$ for which invariant  3-forms exist. The above two only occur as isotropy group.
\end{note*}

\begin{thm}[Main classification theorem]\label{thm_strict_types_and_isotropy_groups}
Let $(M^6,g,J)$ be an almost hermitian structure with parallel torsion. Then the isotropy group of the torsion form is as given in Table \ref{tbl_strict_types_and_isotropy_groups}, and the torsion form belongs to one of two families of forms: 
\begin{gather*}
T^c=\alpha_1\{(e_{14}+e_{23})\wedge e_5\}+(e_{12}-e_{34})\wedge (\alpha_3e_5+\alpha_4e_6)+\alpha_5\{(e_{12}+e_{34})\wedge e_5\},\\
\intertext{or}
T^c=\alpha_1\{(e_{14}+e_{23})\wedge e_5+(e_{13}-e_{24})\wedge e_6\}+\alpha_2\{(-e_{13}+e_{24})\wedge e_5+(e_{14}+e_{23})\wedge e_6\}\\
+\beta_1\{(e_{12}-e_{34})\wedge e_5\}+\beta_2\{(e_{13}-e_{24})\wedge e_5+(e_{14}+e_{23})\wedge e_6\}.
\end{gather*}
Here $e=(e_1,...,e_6)$ is an adapted local frame, and $\alpha_1,\alpha_2,\alpha_3,\alpha_4,\alpha_5,\beta_1,\beta_2$ are real numbers subject to the conditions I--VI and VII--XI of Table \ref{tbl_torsion_form_parameters} in the first respectively second case. The forms determine distinct orbits.
\end{thm}

\noindent
If $T^c$ is as in the theorem, then, in the first case,
\begin{gather*}
\|T_2\|^2=\alpha_1^2,\quad\|T_{12}\|^2=\alpha_1^2+2(\alpha_3^2+\alpha_4^2),\quad\|T_6\|^2=2\alpha_5^2,
\end{gather*}
and, in the second,
\begin{gather*}
\|T_2\|^2=4(\alpha_1^2+\alpha_2^2),\quad\|T_{12}\|^2=2\beta_1^2+4\beta_2^2,\quad\|T_6\|^2=0.
\end{gather*}

\begin{proof}
The theorem is a collection of the results given in Theorem \ref{thm_isotropy_group_W1+W3+W4}, Theorem \ref{thm_singular_orbits_W3}, Corollary \ref{cor_singular_orbits_W1+W3}, Theorem \ref{thm_W3_1dim_isotropy_group}, Theorem \ref{thm_W1W3_singular_orbits_reduced} and Corollary \ref{cor_isotropy_group_W1+W3} below. (For structures of strict type $W_1$ and strict type $W_4$ the result is well known. It follows from our knowledge of the orbit structure of the $U(3)$-modules $\Lambda^3_2$ and $\Lambda^3_6$.) 
\end{proof}

\begin{NB*}
There are now two alternative proofs of the fact that the isotropy group is always non-trivial: a topological proof by Thomas Friedrich, and an algebraic proof by Andrew Swann (both unpublished).
\end{NB*}

\begin{table}
\bdm
\begin{array}{c|c|c}
\hline
\text{strict type} & Iso_o(T^c) & Iso_o(T^c)=Iso(T^c)
\spacedhline
\hline
W_1 & SU(3) & \text{yes} \\
W_3 & U(2)_1,\;SO(3),\;T^2 & \text{no} \\
W_4 & U(2)_0 & \text{yes} \\
W_1\op W_3 & SU(2),\;SO(3),\;T^1 & \text{no} \\
W_3\op W_4 & T^2 & \text{yes}\\
W_1\op W_3\op W_4 & SU(2),\;T^1 & \text{yes}
\spacedhline
\end{array}
\edm
\caption{}
\label{tbl_strict_types_and_isotropy_groups}
\end{table}

\begin{table}
\bdm
\begin{array}{c|c|c|c}
\hline
 & \text{strict type} & Iso_o(T^c) & \text{conditions} 
\spacedhline
\hline
I & W_4 & U(2)_0 & \alpha_1=\alpha_3=\alpha_4=0, \;\alpha_5>0  \\
II & W_1\op W_3 & SU(2) &  \alpha_1>0, \;\alpha_3=\alpha_4=0, \;\alpha_5=0  \\
III & W_1\op W_3 & T^1 &  
\begin{cases} 
\alpha_1>0, \;\alpha_3>0, \;\alpha_4\in\R, \;\alpha_5=0 \\
\alpha_1>0, \;\alpha_3=0, \;\alpha_4>0, \;\alpha_5=0 
\end{cases}\\
IV & W_3\op W_4 & T^2 & 
\begin{cases} 
\alpha_1=0, \;\alpha_3>0, \;\alpha_4\in\R, \;\alpha_5>0  \\
\alpha_1=0, \;\alpha_3=0, \;\alpha_4>0, \;\alpha_5>0
\end{cases}\\
V & W_1\op W_3\op W_4 & SU(2) & \alpha_1>0, \;\alpha_3=\alpha_4=0, \;\alpha_5>0  \\
VI & W_1\op W_3\op W_4 & T^1 & 
\begin{cases} 
\alpha_1>0, \;\alpha_3>0, \;\alpha_4\in\R, \;\alpha_5>0 \\
\alpha_1>0, \;\alpha_3=0, \;\alpha_4>0, \;\alpha_5>0 
\end{cases}
\spacedhline
VII & W_1 & SU(3) & \alpha_1>0,\;\alpha_2=0, \;\beta_1=\beta_2=0  \\
VIII & W_3 & U(2)_1 & \alpha_1=\alpha_2=0, \;\beta_1=0, \;\beta_2\neq0  \\
IX & W_3 & T^2 &  \alpha_1=\alpha_2=0, \;\beta_1\neq0, \;\beta_2=0  \\
X & W_3 & SO(3) &  \alpha_1=\alpha_2=0, \;\beta_1=2\beta_2\neq0 \\
XI & W_1\op W_3 & SO(3) & (\alpha_1,\alpha_2)\neq(0,0), \;\beta_1=2\beta_2\neq0
\spacedhline
\end{array}
\edm
\caption{}
\label{tbl_torsion_form_parameters}
\end{table}

%\newpage
%---------------------------------------------------------------------------- 
\subsection{Non-vanishing divergence}
%---------------------------------------------------------------------------- 
\noindent
Let $(M^6,g,J)$ be an almost hermitian structure with parallel torsion and non-vanishing divergence, $\delta\Omega=2X\neq0$ (cf. Theorem \ref{thm_type_and_diff_eqs}). The vector fields $X$ and $JX$ then are non-zero, parallel and of constant length. The tangent bundle splits orthogonally,
\begin{gather*}
TM^6=\mcH\op\mcV,
\end{gather*}
into two complex distributions: a horizontal distribution $\mcH$ and a vertical distribution $\mcV$, with $\mcV$ spanned by the vector fields $X$ and $JX$. Corresponding to this splitting, almost complex structure and K\"ahler form decompose, $J=J_{\mcH}+J_{\mcV}$ and $\Omega=\Omega_{\mcH}+\Omega_{\mcV}$. We will consider the reduction
\begin{gather*}
\mcR'=\{e=(e_1,...,e_6)\in\mcR:\;e_5=X/\| X\|,\;e_6=JX/\| X\|\}.
\end{gather*}
The structure group is $U(2)$ embedded through $U(2)_0$. 

To begin with, we wish to describe the $U(2)$-modules $\Lambda^3_2$ and $\Lambda^3_{12}$. Let $\R^6=(\R^6,g,J)$ be the standard hermitian space, $\Omega$ the K\"ahler form, and $(e_1,...,e_6)$ an adapted basis. Let $\R^4=(\R^4,g,J_{\mcH})$ and $\R^2=(\R^2,g,J_{\mcV})$ be the hermitian subspaces with basis $(e_1,...,e_4)$ and $(e_5,e_6)$, and K\"ahler forms $\Omega_{\mcH}$ and $\Omega_{\mcV}$ respectively. 

Let us recall the structure of the real $U(2)$-module $\Lambda^2(\R^4)$. The module is reducible,
\begin{gather*}
\Lambda^2(\R^4)=\Lambda^2_+\op\Lambda^2_-=\R\op\m^2\op\Lambda^2_-.
\end{gather*}
The first equality gives the decomposition as $SO(4)$-module into selfdual and anti-selfdual 2-forms, $\Lambda^2_+$ and $\Lambda^2_-$, and the second equality the decomposition as $U(2)$-module. If $\hodge$ denotes the Hodge operator on $\R^4$, then
\begin{gather*}
\Lambda^2_{\pm}=\{\omega\in\Lambda^2(\R^4): \hodge\omega=\pm\omega\},\quad \R=\{\omega\in\Lambda^2(\R^4): \hodge\omega=\omega,J_{\mcH}\omega=\omega\},\\
\m^2=\{\omega\in\Lambda^2(\R^4): \hodge\omega=\omega,J_{\mcH}\omega=-\omega\}.
\end{gather*}
The $U(2)$-module $\m^2$ is isomorphic to $(\C,det^2)$, that is, $a\cdot z=det(a)^2z$ for any $a\in U(2)$, $z\in\C $. The complex structure $F:\m^2\ra\m^2$ is given by
\begin{gather*}
F(\omega)=\frac{1}{2}[\Omega_{\mcH},\omega].
\end{gather*}

We define maps,\label{linalg_Lambda2R4}
\begin{gather*}
i_1:\m^2\ra\Lambda^3_2,\quad\omega\mapsto\frac{1}{\| X\|}\{\omega\wedge X-F\omega\wedge JX\},\\
\intertext{and}
i_2:\m^2\ra\Lambda^3_{12},\quad\omega\mapsto\frac{1}{\| X\|}\{\omega\wedge X+F\omega\wedge JX\},\\
i_3:\Lambda^2_-\ox\C\ra\Lambda^3_{12},\quad\omega_1+i\omega_2\mapsto\frac{ 1}{\| X\|}\{\omega_1\wedge X+\omega_2\wedge JX\},\\
i_5:\R^4\ra\Lambda^3_{12},\quad V\mapsto (\Omega_{\mcH}-\Omega_{\mcV})\wedge V.
\end{gather*}

\begin{prop}\label{prop_Lambda3s12_U20_decomposition}
As $U(2)$-modules,
\begin{gather*}
\Lambda^3_2=\m^2,\quad\Lambda^3_{12}=\m^2\op(\Lambda^2_-\ox\C)\op\R^4.
\end{gather*}
The $U(2)$-modules $\m^2$ and $\R^4$ are real irreducible, the $U(2)$-module $\Lambda^2_-\ox\C$ is complex irreducible.
\end{prop}

\begin{proof}
The maps $i_1$, $i_2$, $i_3$ and $i_5$ are well defined (i.e. take values in the respective spaces), injective, $U(2)$-equivariant, and their images have trivial intersection.
\end{proof}

\noindent
The orbit structure of $\m^2$ and $\R^4$ is obvious, that of $\Lambda^2_-\ox\C$ may be seen easily: $U(2)$ acts on $\Lambda^2_-=\su(2)$ through the adjoint representation, the kernel of this action equals $SO(2)$, its image equals $SO(3)$. As $SO(3)$-module, $\Lambda^2_-\ox\C$ is isomorphic to $\R^3\ox\C$, the complexification of the defining representation of $SO(3)$. The orbit structure is as follows.

\begin{prop}\label{prop_orbit_structure_SO3_R6}
The $SO(3)$-module $\R^6=\R^3\op\R^3$ contains precisely three kinds of orbits: the zero orbit with isotropy group $SO(3)$, orbits with isotropy group $SO(2)$, and orbits with trivial isotropy group. If $(e_1,e_2,e_3)$ is an orthonormal basis in $\R^3$, orbits with isotropy group $SO(2)$ may be parameterized by pairs of vectors $(\lambda e_1,\mu_1 e_1)$ as determined by the two sets of parameters $\lambda>0$, $\mu_1\in\R$ and $\lambda=0$, $\mu_1>0$; orbits with trivial isotropy group may be parameterized by pairs of vectors $(\lambda e_1,\mu_1 e_1+\mu_2 e_2)$ as determined by the set of parameters $\lambda>0$, $\mu_1\in\R$, $\mu_2>0$.
\end{prop}

The decomposition of $\Lambda^3_2$ and $\Lambda^3_{12}$ under the $U(2)$-action determines a splitting of the corresponding bundles. This allows us to express the 3-forms $T_2$ and $T_{12}$ in terms of suitable (horizontal) 2-forms. Let 
\begin{gather*}
T_2=i_1(\Omega_1),\quad T_{12}=i_2(\Omega_2)+i_3(\Omega_3+i \Omega_4)+i_5(Y)
\end{gather*}
for forms $\Omega_1,\Omega_2,\Omega_3,\Omega_4$ and $Y$.

\begin{lem}\label{lem_restriction_classes}
$\Omega_1=\Omega_2$ and $Y=0$.
\end{lem}

\begin{proof}
The torsion form is parallel and thus coclosed. It follows that $\delta^g(T_2-T_{12})=0$, \cite[Proposition 4.1]{AlexandrovFriedrichSchoemann:2005}. Clearly $\delta^c(T_2-T_{12})=0$. For a $k$-form $\omega$, $k\geq2$, the difference of the two codifferentials is given by
\begin{gather*}
(\delta^g-\delta^c)\omega=\frac{1}{2}\sum_{i,j=1}^6(e_{ij}\haken T^c)\wedge (e_{ij}\haken\omega),
\end{gather*}
\cite[Proposition 5.1]{AgricolaFriedrich:2004}. Letting $\omega=T_2-T_{12}$, the left-hand side vanishes. The resulting 15 equations give the result.
\end{proof}

\begin{NB}
In the case of vanishing divergence, $\delta\Omega=0$, the equations in the preceding proof are always satisfied.
\end{NB}

A consequence of Lemma \ref{lem_restriction_classes} is the following.

\begin{thm}[{\cite[Theorem 4.2]{AlexandrovFriedrichSchoemann:2005}}]\label{thm_non_existance_W1+W4}
Almost hermitian structures with parallel torsion of strict type $W_1\op W_4$ do not exist.
\end{thm}

\begin{proof}
Suppose there exists an almost hermitian structure with parallel torsion of strict type $W_1 \op W_4$. This structure has non-vanishing divergence, by definition. Now $T_{12}=0$ and Lemma \ref{lem_restriction_classes} implies $T_2=0$, a contradiction.
\end{proof}

The following is the central result of this section.

\begin{thm}\label{thm_isotropy_group_W1+W3+W4}
Let $(M^6,g,J)$ be an almost hermitian structure with parallel torsion and non-vanishing divergence. Then there exists a local frame such that
\begin{gather*}
T^c=\alpha_1\{(e_{14}+e_{23})\wedge e_5\}+(e_{12}-e_{34})\wedge (\alpha_3e_5+\alpha_4e_6)+\alpha_5\{(e_{12}+e_{34})\wedge e_5\},
\end{gather*}
with $\alpha_1,\alpha_3,\alpha_4,\alpha_5$ real numbers subject to the conditions I, IV, V, VI of Table \ref{tbl_torsion_form_parameters}. The forms determine distinct orbits, and the isotropy group is as given in Table \ref{tbl_strict_types_and_isotropy_groups}.
\end{thm}

\begin{proof}
From our knowledge of the $U(2)$-orbit structure of the modules $\m^2$ and $\Lambda^2_-\ox\C$, it follows that we may find a frame $e\in\mcR'$ such that
\begin{gather*}
\Omega_1=\alpha_1(e_{14}+e_{23}),\quad\Omega_3+i\Omega_4=(\alpha_3+i\alpha_4)(e_{12}-e_{34})+i\gamma(e_{13}+e_{24}),
\end{gather*}
with $\alpha_1,\alpha_3,\alpha_4,\gamma$ real numbers. Moreover, we may choose $\alpha_1$ non-negative, and restrict the parameters $\alpha_3$, $\alpha_4$ and $\gamma$ according to the description of the $SO(3)$-orbits of $\R^6=\R^3\op\R^3$ given in Proposition \ref{prop_orbit_structure_SO3_R6}, $(\alpha_3,\alpha_4,\gamma)\widehat{=}(\lambda,\mu_1,\mu_2)$. 

For $\gamma=0$, compute $T_2$ and $T_{12}$ using the maps $i_1$, $i_2$ and $i_3$, and notice that 
\begin{gather*}
T_6=\Omega\wedge X=\alpha_5(e_{12}+e_{34})\wedge e_5
\end{gather*}
with $\alpha_5=\|X\|=\frac{1}{\sqrt{2}}\|T_6\|$ a positive real number. This yields $T^c=T_2+T_{12}+T_6$ with isotropy groups and parameters as stated. 

For $\gamma\neq 0$, notice that $Iso(\Omega_3+i\Omega_4)\cap SO(3)=T^2\cap SO(3)\{e\}$, where $SO(3)\subset U(2)=T^1\x SO(3)$. The holonomy group, therefore, must be either trivial or $T^1$. In both cases the Bianchi identities imply $\gamma=0$, and the case $\gamma\neq 0$, therefore, cannot occur.
\end{proof}

%\newpage
%---------------------------------------------------------------------------- 
\subsection{Vanishing divergence I: Singular orbits in $\Lambda^3_{12}$ and $\Lambda^3_2\op\Lambda^3_{12}$}
%---------------------------------------------------------------------------- 
An orbit is called \emph{singular} if its isotropy group is non trivial, that is, if the dimension $dim\; Iso(T^c)\geq1$. An orbit is called \emph{$T^1$-singular} if $T^1\subset Iso(T^c)$, where $T^1$ is a torus in a given conjugacy class of tori. Every singular orbit is $T^1$-singular for some $T^1$.

Fix a conjugacy class of tori of dimension one, and  let $T^1$ be a torus in that class. Any $T^1$-singular orbit intersects the space $V=(\Lambda^3_{12})^{T^1}$, the trivial component of the $T^1$-module $\Lambda^3_{12}$. The invariance group is
\begin{gather*}
Inv(V)=\{a\in U(3): aV\subset V\}.
\end{gather*}
The problem of describing the $T^1$-singular orbits in $\Lambda^3_{12}$ is equivalent to describing the $T^1$-singular $Inv(V)$-orbits in $V$. While the invariance group is difficult to determine in general, finding a large subgroup $G\subset Inv(V)$ is easier. A description of the orbit structure of the $G$-module $V$ yields an ``upper bound'' for the $T^1$-singular orbits in $\Lambda^3_{12}$, that is, a description of orbits that may not be unique but from which no orbit is missing. (Restricting to the singular, or $T^1$-singular orbits, a sharper bound may be obtained.)

There are infinitely many conjugacy classes of 1-dimensional tori. It is possible, however, to choose an infinite set of tori $T^1$, one from each conjugacy class, such that only finitely many spaces $V=(\Lambda^3_{12})^{T^1}$ occur. A description of the orbit structure of the various $V$ then leads to an upper bound for the singular orbits in $\Lambda^3_{12}$. This is an outline of our approach to obtain a description of the singular orbits in $\Lambda^3_{12}$ and $\Lambda^3_2\op\Lambda^3_{12}$, details follow.

To begin with, we describe the conjugacy classes of tori of dimensional one in $U(3)$.  Let 
\begin{gather*}
T^1_{(k_1,k_2,k_3)}=\{diag(e^{i k_1 t},e^{i k_2 t},e^{i k_3 t}):t\in\R\}
\end{gather*}
 be the torus in $U(3)$ defined by a tuple of integers $(k_1,k_2,k_3)\in\Z^3\backslash 0$, and let $\bar{\Delta}=(\Z^3\backslash 0)/\sim$ denote the quotient in $\Z^3\backslash 0$ where tuples $(k_1,k_2,k_3)$ and $(k'_1,k'_2,k'_3)$ are identified if $(k'_1,k'_2,k'_3)=\lambda(k_{\sigma(1)},k_{\sigma(2)},k_{\sigma(3)})$ for some integer $\lambda\in\Z\backslash 0$ and permutation $\sigma\in S_3$, the permutation group of three letters.

\begin{lem}
The $U(3)$-conjugacy classes of 1-dimensional tori in $U(3)$ are in one-to-one correspondence with elements of $\bar{\Delta}$.
\end{lem}

We now choose a unique element from each conjugacy class. The number of zeros in tuples representing a given element in $\bar{\Delta}$ is independent of the choice of the tuple. Therefore, if the number of zeros is two, one may find tuples $(k_1,k_2,k_3)$ such that $k_1>0$ and  $k_2=k_3=0$, if the number of zeros is one, such that $k_1>0$, $k_2\neq0$ and  $k_3=0$, and, if the number of zeros is zero, such that $k_1\geq k_2>0$ and  $k_2\geq k_3$. Inside the set of tuples with this ordering there is a unique tuple $(k_1,k_2,k_3)$ such that the greatest common divisor of $k_1$, $k_2$ and $k_3$ is one. Here, any natural number $d\in\N$ divides zero, and $d$ divides $-k$, $k>0$, if and only if $d$ divides $k$. The set of tuples thus singled out we denote by $\Delta$.

Now let
\begin{gather*}
\begin{aligned}
\Delta_1&=\{(k_1,0,0)\in\Delta:\;k_1>0\},&\Delta_2&=\{(k_1,k_2,0)\in\Delta:\;k_1>0,\;k_2=k_1\},\\
\Delta_3&=\{(k_1,k_2,0)\in\Delta:\;k_1>0,\;k_2=-k_1\},&\Delta_4&=\{(k_1,k_2,0)\in\Delta:\;k_1>0,\;k_2\neq\pm k_1\},
\end{aligned}\\
\begin{aligned}
\Delta_5&=\{(k_1,k_2,k_3)\in\Delta:\;k_1>k_2>0,\;k_2>k_3,\;k_3=-(k_1-k_2)\},\\
\Delta_6&=\{(k_1,k_2,k_3)\in\Delta:\;k_1>k_2>0,\;k_2>k_3,\;k_3=+(k_1-k_2)\},\\
\Delta_7&=\{(k_1,k_2,k_3)\in\Delta:\;k_1\geq k_2>0,\;k_2\geq k_3,\;k_3\neq\pm(k_1-k_2),\;k_3=-(k_1+k_2)\},\\
\Delta_8&=\{(k_1,k_2,k_3)\in\Delta:\;k_1\geq k_2>0,\;k_2\geq k_3,\;k_3\neq\pm(k_1-k_2),\;k_3\neq-(k_1+k_2)\}.
\end{aligned}
\end{gather*}

\begin{lem}
The sets $\Delta_i$ are pairwise disjoint and $\Delta=\Delta_1\cup...\cup\Delta_8$.
\end{lem}

We write $T^1=T^1_{(k_1,k_2,k_3)}\in\Delta_i$ for a torus determined by integers $(k_1,k_2,k_3)\in\Delta_i$. 

\begin{lem}\label{lem_properties_families_of_tori}
For each $i=1,...,8$, the spaces $(\Lambda^3_{2})^{T^1}$, $(\Lambda^3_{12})^{T^1}$ and $(\Lambda^3_6)^{T^1}$ are independent of the choice of $T^1\in\Delta_i$, .
\end{lem}

\noindent
The partition of $\Delta$ was motivated by Lemma \ref{lem_properties_families_of_tori}. Notice that the $\Delta_i$ are distinguishable by the dimensions of the spaces $(\Lambda^3_{2})^{T^1}$, $(\Lambda^3_{12})^{T^1}$ and $(\Lambda^3_6)^{T^1}$, with the exception of $\Delta_5$ and $\Delta_6$, see Table \ref{tbl_trivial_T1_submodules}. 

\begin{table}
\bdm
\begin{array}{r||cccccccc}
\hline
& \Delta_1 & \Delta_2 & \Delta_3 & \Delta_4 & \Delta_5 & \Delta_6 & \Delta_7 & \Delta_8 
\spacedhline
\hline
dim\;(\Lambda^3_{2})^{T^1} & 0 & 0 & 2 & 0 & 0 & 0 & 2 & 0 \\
dim\;(\Lambda^3_{12})^{T^1} & 4 & 6 & 4 & 2 & 2 & 2 & 0 & 0 \\
dim\;(\Lambda^3_{6})^{T^1} & 4 & 2 & 2 & 2 & 0 & 0 & 0 & 0 
\spacedhline
\end{array}
\edm
\caption{}
\label{tbl_trivial_T1_submodules}
\end{table}

The next step is to understand the orbit structure of the $G$-modules $(\Lambda^3_{12})^{T^1}$ with $T^1\in\Delta_1,...,\Delta_8$, and $G$ either the invariance group of that module or some large subgroup of it. The space $(\Lambda^3_{12})^{T^1}$ with $T^1\in\Delta_i$ is non-trivial for $\Delta_1,...,\Delta_6$, and we find that the invariance group contains a group $G$ as follows: 
\begin{gather*}
\Delta_1: G=U(2),\quad\Delta_2: G=T^1\x SU(2),\quad\Delta_3: G=T^1\x T^1,\quad\Delta_4,\Delta_5,\Delta_6: G=T^1. 
\end{gather*}
The representation of $G$ on $(\Lambda^3_{12})^{T^1}$ can be determined (all groups are well known and the representations are of low dimension). The group $SU(2)$, for instance, acts on the 6-dimensional $(\Lambda^3_{12})^{T^1}$ through the complexification of the (real) adjoint representation. We thus find the orbit structure of the six $G$-modules $(\Lambda^3_{12})^{T^1}$. The orbits are described by families of forms as follows:
\begin{align*}
&\Delta_1,\Delta_2,\Delta_4: &T_{12}=&\;\beta_1\{(e_{12}-e_{34})\wedge e_5\}, \\
&\Delta_3: &T_{12}=&\;\beta_1\{(e_{12}-e_{34})\wedge e_5\} \\
& &&\;+\beta_2\{(e_{13}-e_{24})\wedge e_5+(e_{14}+e_{23})\wedge e_6\}, \\
&\Delta_5: &T_{12}=&\;\beta_1\{(e_{13}+e_{24})\wedge e_5+(e_{14}-e_{23})\wedge e_6\}, \\
&\Delta_6: &T_{12}=&\;\beta_1\{(e_{13}+e_{24})\wedge e_5-(e_{14}-e_{23})\wedge e_6\}.
\end{align*}
Here $\beta_1,\beta_2$ are non-negative real numbers.

The final step is to remove redundant forms. The isotropy groups of the forms in case $\Delta_3$ are as stated in Theorem \ref{thm_singular_orbits_W3} below, which is easy to verify. The isotropy groups of forms in case $\Delta_4$ and $\Delta_5$ is $U(2)_1$. The $U(3)$-orbits of forms in $\Lambda^3_{12}$ whose isotropy group equals $U(2)_1$ are parameterized by one positive real number, a fact not difficult to check separately. The orbits defined by the forms in case $\Delta_5$ and $\Delta_6$, therefore, are contained in the family of orbits defined in case $\Delta_3$. 

\begin{thm}\label{thm_singular_orbits_W3}
Apart from the trivial orbit, the singular orbits of the $U(3)$-module $\Lambda^3_{12}$ are contained in the family of orbits determined by
\begin{gather*}
T_{12}=\beta_1\{(e_{12}-e_{34})\wedge e_5\}+\beta_2\{(e_{13}-e_{24})\wedge e_5+(e_{14}+e_{23})\wedge e_6\}
\end{gather*}
with $\beta_1,\beta_2$ real numbers, both non-negative, not both zero. The isotropy group is as follows: (i) $\beta_1\neq0$, $\beta_2=0$: $Iso_o(T_{12})=T^2$, (ii) $\beta_1=0$, $\beta_2\neq0$: $Iso_o(T_{12})=U(2)_1$, (iii) $\beta_1\neq0$, $\beta_2\neq0$, $\beta_1=2\beta_2$: $Iso_o(T_{12})=SO(3)$, (iv) $\beta_1\neq0$, $\beta_2\neq0$, $\beta_1\neq2\beta_2$: $Iso_o(T_{12})=T^1$. The forms determine distinct orbits in case (i),(ii),(iii), in case (iv) a small ambiguity remains.
\end{thm}

\noindent
The ambiguity of case (iv) is explained in Remark \ref{rmk_Z_2_ambiguity} below.

\begin{cor}\label{cor_singular_orbits_W1+W3}
The singular orbits of the $U(3)$-module $\Lambda^3_2\op\Lambda^3_{12}$ that intersect $\Lambda^3_2$ and $\Lambda^3_{12}$ are contained in the family determined by $T^c=T_2+T_{12}$ with $T_2$ a general form in $\Lambda^3_2$ and $T_{12}$ a form in $\Lambda^3_{12}$ as in Theorem \ref{thm_singular_orbits_W3}. That is,
\begin{gather*}
T_2=\;\alpha_1\{(e_{14}+e_{23})\wedge e_5+(e_{13}-e_{24})\wedge e_6\}+\alpha_2\{(-e_{13}+e_{24})\wedge e_5+(e_{14}+e_{23})\wedge e_6\},\\
T_{12}=\;\beta_1\{(e_{12}-e_{34})\wedge e_5\}+\beta_2\{(e_{13}-e_{24})\wedge e_5+(e_{14}+e_{23})\wedge e_6\},
\end{gather*}
with $\alpha_1,\alpha_2,\beta_1,\beta_2$ real numbers, $\alpha_1,\alpha_2$ not both zero, $\beta_1,\beta_2$ both non-negative, not both zero. The isotropy group is as follows: (i) $\beta_1\neq0$, $\beta_2=0$: $Iso_o(T^c)=T^1$, (ii) $\beta_1=0$, $\beta_2\neq0$: $Iso_o(T^c)=SU(2)$, (iii) $\beta_1\neq0$, $\beta_2\neq0$, $\beta_1=2\beta_2$: $Iso_o(T^c)=SO(3)$, (iv) $\beta_1\neq0$, $\beta_2\neq0$, $\beta_1\neq2\beta_2$: $Iso_o(T^c)=T^1$. In case (i),(ii) restriction to parameters with $\alpha_1>0$, $\alpha_2=0$ yields a complete set of forms which determine distinct orbits. In case (iii) the forms determine distinct orbits, in case (iv) a small ambiguity remains.
\end{cor}

\noindent
In case (i),(ii) the isotropy group of $T_{12}$ is either $U(2)_1$ or $T^2$. It acts on $\Lambda^3_2$ through rotation. This shows, that a complete set of distinct orbits is determined by restricting parameters as stated. The case (iii),(iv) is explained in the following remark.

\begin{NB}\label{rmk_Z_2_ambiguity}
Invariant polynomials may be used to show that a set of forms determines distinct orbits. Since the exponential map $\u(3)\ra U(3)$ is surjective, the problem of finding $U(3)$-invariant polynomials on any finite dimensional $U(3)$-module can be linearized. For $\Lambda^3_2\op\Lambda^3_{12}$, using a computer, one may find that the space of invariant polynomials of degree $\leq4$ is of dimension 8. Making use of such polynomials, one can show that in Corollary \ref{cor_singular_orbits_W1+W3}, case (iii), the forms determine distinct orbits, and that in Theorem \ref{thm_singular_orbits_W3}, case (iv), and Corollary \ref{cor_singular_orbits_W1+W3}, case (iv), precisely two forms may determine the same orbit, while all other forms determine distinct orbits.
\end{NB}

%\newpage
%---------------------------------------------------------------------------- 
\subsection{Vanishing divergence II: Further restrictions}
%---------------------------------------------------------------------------- 
The first Bianchi identities further restrict the set of orbits that may occur as the orbit of the torsion form (see page \pageref{page_Biachi_identities}).

\begin{thm}[{\cite[Theorem 4.8]{AlexandrovFriedrichSchoemann:2005}}]\label{thm_W3_1dim_isotropy_group}
Almost hermitian structures with parallel torsion of strict type $W_3$ and $dim\;Iso(T^c)=1$ do not exist.
\end{thm}

\begin{proof}
Suppose that $dim\;Iso(T^c)=1$. Then $Iso_o(T^c)=T^1$, by Theorem \ref{thm_singular_orbits_W3}. If $e\in\mcR$ is a frame as in that theorem, then
\begin{align*}
\sigma_{T^c}=(-\beta_1^2+2\beta_2^2)\;e_{1234}-2\beta_2^2\;(e_{1256}+e_{3456}). 
\end{align*}
On the other hand, $\hol(\nabla^c)\subset\iso(T^c)=Lin_{\R}\{e_{12}-e_{34}\}$, and therefore
\begin{align*}
(\cycsum R^c)(e_1,e_2,e_5,e_6)=0. 
\end{align*}
The Bianchi identities, $\cycsum R^c=\sigma_{T^c}$, imply $\beta_2=0$. Now $Iso_o(T^c)=T^2$, by Theorem \ref{thm_singular_orbits_W3}, a contradiction.
\end{proof}

\begin{thm}\label{thm_W1W3_singular_orbits_reduced}
Let $(M^6,g,J)$ be an almost hermitian structure with parallel torsion of strict type $W_1\op W_3$ and $Iso_o(T^c)=SU(2)$ or $Iso_o(T^c)=T^1$. Then there exists a local frame such that
\begin{align*}
T^c=\alpha_1\{(e_{14}+e_{23})\wedge e_5\}+(e_{12}-e_{34})\wedge(\alpha_3e_5+\alpha_4e_6)
\end{align*}
with $\alpha_1,\alpha_3,\alpha_4$  real numbers subject to the conditions II, III of Table \ref{tbl_torsion_form_parameters}. The forms determine distinct orbits, and the isotropy group is as given in Table \ref{tbl_strict_types_and_isotropy_groups}.
\end{thm}

\begin{proof}
Let $e\in\mcR$ be a frame as in Corollary \ref{cor_singular_orbits_W1+W3}. The torsion form $T^c$ is contained in the 6-dimensional space of 3-forms $T=T_2+T_{12}$, $T_{12}=T'_{12}+T''_{12}$ with 
\begin{gather*}
T_2=\alpha_1\{(e_{14}+e_{23})\wedge e_5+(e_{13}-e_{24})\wedge e_6\}+\alpha_2\{(-e_{13}+e_{24})\wedge e_5+(e_{14}+e_{23})\wedge e_6\},\\
T'_{12}=\beta_1\{(e_{14}+e_{23})\wedge e_5-(e_{13}-e_{24})\wedge e_6\}+\beta_2\{(-e_{13}+e_{24})\wedge e_5-(e_{14}+e_{23})\wedge e_6\},\\
T''_{12}=(e_{12}-e_{34})\wedge(\alpha_3e_5+\alpha_4e_6),
\end{gather*}
where $\alpha_1,\alpha_2,\beta_1,\beta_2,\alpha_3,\alpha_4$ are real numbers. Let us determine a new 3-form in the orbit determined by the torsion form. The invariance group of the 2-dimensional space of forms $T''_{12}$ contains $T^2$. Now $T^2$ acts invariantly on the space of forms $T'_{12}$ such that we may assume $\alpha_1>0$, $\alpha_2=0$, $\beta_1\geq0$ and $\beta_2=0$. A calculation gives 
\begin{gather*}
\sigma_{T^c}=(2\alpha_1^2+2\beta_1^2-\alpha_3^2-\alpha_4^2)\;e_{1234}+2(\alpha_1^2-\beta_1^2)(e_{1256}+e_{3456}),
\intertext{and obviously}
(\cycsum R^c)(e_1,e_2,e_5,e_6)=0. 
\end{gather*}
The Bianchi identities, $\cycsum R^c=\sigma_{T^c}$, imply $\alpha_1^2-\beta_1^2=0$, and therefore $\alpha_1=\beta_1$. The torsion form $T^c$ is thus as stated.

Now 
\begin{align*}
Iso_o((e_{14}+e_{23})\wedge e_5)=&\;Iso_o((e_{14}+e_{23})\wedge e_5+(e_{13}-e_{24}\wedge e_6))\\
&\;\cap Iso_o((e_{14}+e_{23})\wedge e_5-(e_{13}-e_{24})\wedge e_6)\\
=&\;SU(3)\cap U(2)_1=SU(2), 
\end{align*}
and $T^2\subset Iso((e_{12}-e_{34})\wedge(\alpha_3e_5+\alpha_4e_6))$ for any $\alpha_3,\alpha_4\in\R$. Therefore, $Iso_o(T^c)=SU(2)$ if and only if $\alpha_3=\alpha_4=0$ (case II). The group $SU(2)$ does not act invariantly on the space of forms $T''_{12}$; it does, however, act invariantly on a 6-dimensional space which contains that space. The action is through the adjoint representation. The remaining part of the theorem now follows from our description of the $SO(3)$-orbits in $\R^6=\R^3\op\R^3$, see Proposition \ref{prop_orbit_structure_SO3_R6} (case III).
\end{proof}

%\newpage
%---------------------------------------------------------------------------- 
\subsection{Vanishing divergence III: Trivial isotropy groups}
%---------------------------------------------------------------------------- 
\noindent
The module $\Lambda^3_2\op\Lambda^3_{12}$ contains many orbits with trivial isotropy group ($dim\;Iso(T^c)=0$). These orbits do not occur as the orbit of the torsion form.

\begin{prop}
Let $(M^6,g,J)$ be an almost hermitian structure with parallel torsion of type $W_1\op W_3$ and trivial holonomy. Then $dim\;Iso(T^c)\geq1$, and $3\|T_2\|^2=\|T_{12}\|^2$. 
\end{prop}

\begin{proof}
Trivial isotropy group implies trivial holonomy. The Bianchi identities, in this situation, are algebraic relations in the coefficients of the torsion form. Encoded in the Clifford algebra, the Bianchi identities hold if and only if the square of the torsion form inside the Clifford algebra is a scalar, $(T^c)^2\in\R\subset Cl(\R^6)$, (cf. Remark \ref{rmk_Sternberg_criteria}). We will choose a family of 3-forms such that any orbit in $\Lambda^3_2\op\Lambda^3_{12}$ is represented, not uniquely in general, and apply the criteria. 

To this end, let $T_2\in\Lambda^3_2$ be a generic 3-form, 
\begin{gather*}
T_2=\alpha_1\{(e_{14}+e_{23})\wedge e_5+(e_{13}-e_{24})\wedge e_6\}+\alpha_2\{(-e_{13}+e_{24})\wedge e_5+(e_{14}+e_{23})\wedge e_6\},
\end{gather*}
with $\alpha_1,\alpha_2$ real numbers. Consider $\Lambda^3_{12}$ as $U(2)$-module. As such it decomposes, $\Lambda^3_{12}=\m^2\op(\Lambda^2_-\ox\C)\op\C^2$, (cf. page \pageref{linalg_Lambda2R4}f), and a generic 3-form $T_{12}\in\Lambda^3_{12}$ may be written as $T_{12}=i_2(\Omega_2)+i_3(\Omega_3+i \Omega_4)+i_5(Y)$ with $\Omega_2,\Omega_3,\Omega_4$ and $Y$ forms of appropriate type. We may choose $\alpha_1\geq0$, $\alpha_2=0$. Since $SU(3)\cap U(2)=SU(2)$, we may further choose $Y$ to depend on one (non-negative) parameter.  A general form $T^c=T_2+T_{12}$ may thus be chosen to depend on 10, rather than 14, real parameters. 

Now $(T^c)^2\in Cl(\R^6)$ must be a scalar. If $Y\neq0$, the corresponding equations show that $\Omega_2$ is determined by $\Omega_3+i \Omega_4$, and that two parameters in $\Omega_3+i \Omega_4$ are zero. A direct computation of the isotropy algebra of $T_{12}$ then shows $dim\; Iso(T_{12})>0$. If $Y=0$, using our knowledge of the $SO(3)$-orbit structure of $\R^6=\R^3\op\R^3$ (Proposition \ref{prop_orbit_structure_SO3_R6}), we may choose three parameters in $\Omega_3+i \Omega_4$ to be zero. In this case then, a general form $T^c=T_2+T_{12}$ depends on 6 real parameters. Again $(T^c)^2\in Cl(\R^6)$ must be a scalar, and we find that that case $Iso(\Omega_3+i \Omega_4)\cap SO(3)=\{e\}$ is not possible. If $Iso(\Omega_3+i \Omega_4)\cap SO(3)\neq\{e\}$ then $dim\; Iso(T_{12})>0$. Hence, if $T^c=T_2+T_{12}$ is a torsion form of an almost hermitian structure with parallel torsion, the orbit generated by $T_{12}$ is singular.

The singular orbits in $\Lambda^3_{12}$ are parameterized by two non-negative real numbers, Theorem \ref{thm_singular_orbits_W3}. We replace $T_{12}$ in $T^c=T_2+T_{12}$ by the form given there, and $T_2$ by the form given above (two parameters). It is now easy to verify that $(T^c)^2\in Cl(\R^6)$ is a scalar if and only if $T^c$ satisfies the conditions of the proposition.
\end{proof}

\begin{cor}[cf. {\cite[Theorem 4.7]{AlexandrovFriedrichSchoemann:2005}}]
Almost hermitian structures with parallel torsion and discrete holonomy of strict type $W_1$ or strict type $W_3$ do not exist.
\end{cor}

\begin{cor}\label{cor_isotropy_group_W1+W3}
Almost hermitian structures with parallel torsion of type $W_1\op W_3$ such that the isotropy group of the torsion form is trivial do not exist.
\end{cor}

%\newpage
%%%%%%%%%%%%%%%%%%%%%%%%%%%%%%%%%%%%%%%%%%%%%%%%%%%%%%%%%%%%%%%%%%%%%%%%%%%%%
\section{Geometry}
%%%%%%%%%%%%%%%%%%%%%%%%%%%%%%%%%%%%%%%%%%%%%%%%%%%%%%%%%%%%%%%%%%%%%%%%%%%%%
%
\noindent
The classification of the torsion orbits provides a good basis for a systematic investigation into the geometry of almost hermitian structures with parallel torsion. In the following, we treat structures related to the groups $SU(2)$, $SO(3)$, $T^2$ and $T^1$, four major classes.
  
The four classes arrange a large part of all possible geometries in a coarse way. Each class contains structures of various strict type, and with various isotropy and holonomy group. The geometry of structures in a given class admits, to some extent, a uniform description. 

The classes have some intersection with other classes often studied. For the most part, however, structures in these four classes have not been studied before. The results here together with previous results give a very complete picture of almost hermitian structures with parallel torsion. (The problem of describing nearly K\"ahler structures with holonomy $SU(3)$, and generalized Hopf structures with holonomy $U(2)_0$, we did not pursue.)

%\newpage
%----------------------------------------------------------------------------
\subsection{Structures related to $SU(2)$}\label{sec_geometry_SU2}
%----------------------------------------------------------------------------
\noindent
Let $(M^6,g,J)$ be an almost hermitian structure with parallel torsion such that
\begin{gather*}
Hol(\nabla^c)\subset SU(2),\quad SU(2)\subset Iso(T^c).
\end{gather*}
The strict type then is either (i) $W_1\op W_3$, (ii) $W_1\op W_3\op W_4$, or (iii) $W_4$. There exists a local frame such that torsion form 
\begin{gather*}
T^c=\alpha_1\{(e_{14}+e_{23})\wedge e_5\}+\alpha_5\{(e_{12}+e_{34})\wedge e_5\},
\end{gather*}
with $\alpha_1,\alpha_5$ real numbers satisfying one of the following conditions: (i) $\alpha_1>0$, $\alpha_5=0$, (ii)  $\alpha_1>0$, $\alpha_5>0$, (iii)  $\alpha_1=0$, $\alpha_5>0$. Then
\begin{gather*}
\|T_2\|^2=\alpha_1^2=\|T_{12}\|^2,\quad\|T_6\|^2=2\alpha_5^2.
\end{gather*}
(Cf. Theorem \ref{thm_strict_types_and_isotropy_groups}, and cases II, V, I of Table \ref{tbl_torsion_form_parameters}.) Throughout such a local frame is fixed.

The tangent bundle splits orthogonally into a horizontal distribution $\mcH$ and a vertical distribution $\mcV$,
\begin{gather*}
TM^6=\mcH\op\mcV.
\end{gather*}
Locally, $\mcH$ is spanned by the vector fields $e_1,...,e_4$, and $\mcV$ by the vector fields $e_5,e_6$. The distribution $\mcH$ is complex 2-dimensional, the distribution $\mcV$ is complex 1-dimensional, and both are parallel. Let $J=J_{\mcH}+J_{\mcV}$ and $\Omega=\Omega_{\mcH}+\Omega_{\mcV}$ be the corresponding decompositions of almost complex structure and K\"ahler form. 

There exists a associated $S^2$-family of parallel orthogonal almost complex structures. The bundle $\Lambda^2(\mcH)$ of 2-forms over $\mcH$ splits in correspondence to the decomposition of the $SU(2)$-module $\Lambda^2(\R^4)$ into (irreducible) submodules, $\Lambda^2(\R^4)=\Lambda^2_+\op\Lambda^2_-=\R^1\op\m^2\op\Lambda^2_-$, (cf. page \pageref{linalg_Lambda2R4}). The form $\Omega_{\mcH}$ spans the rank one subbundle of $\Lambda^2(\mcH)$ corresponding to $\R^1$, and we may chose two orthogonal parallel forms $\Omega_1,\Omega_2$ that span the subbundle corresponding to $\m^2$. We may assume $\|\Omega_1\|=\|\Omega_2\|=\|\Omega_{\mcH}\|$. For $(q_1,q_2,q_3)\in\R^3$ with $q_1^2+q_2^2+q_3^2=1$, the form $\tilde{\Omega}_{\mcH}=q_1\Omega_1+q_2\Omega_2+q_3\Omega_{\mcH}$ determines an orthogonal endomorphism $\tilde{J}_{\mcH}$ with $(\tilde{J}_{\mcH})^2=-id_{\mcH}$. Let $\tilde{J}=\tilde{J}_{\mcH}+J_{\mcV}$. It follows that $\tilde{J}$ is $g$-orthogonal, $\tilde{J}^2=-id$, and $\nabla^c\tilde{J}=0$. 

The almost hermitian structure $(g,\tilde{J})$ is of type $W_1\op W_3\op W_4$, and $\nabla^c$ is its characteristic connection. The torsion form is parallel, $Hol(\nabla^c)\subset SU(2)$, and $SU(2)\subset Iso(T^c)$. The structure is related to $SU(2)$.

Any orbit can be realized by replacing $J$ through $\tilde{J}$. To see this, let $\tilde{\alpha}_1,\tilde{\alpha}_5$ be another set of parameters. We may assume, $\alpha_1^2+\alpha_5^2=1$ and $\tilde{\alpha}_1^2+\tilde{\alpha}_5^2=1$. With $\tilde{\Omega}_{\mcH}=(\alpha_1\tilde{\alpha}_5-\alpha_5\tilde{\alpha}_1)\Omega_1+(\alpha_1\tilde{\alpha}_1+\alpha_5\tilde{\alpha}_5)\Omega_{\mcH}$, we obtain a family of structures $(g,\tilde{J})$. A $(g,\tilde{J})$-adapted frame is $\tilde{e}_1=-\tilde{J}\tilde{e}_2$, $\tilde{e}_2=e_2$, $\tilde{e}_3=-\tilde{J}\tilde{e}_4$, $\tilde{e}_4=e_4$, $\tilde{e}_5=e_5$, $\tilde{e}_6=e_6$, and with respect to its dual frame 
\begin{align*}
T^c=\tilde{\alpha}_1\{(\tilde{e}_{14}+\tilde{e}_{23})\wedge \tilde{e}_5\}+\tilde{\alpha}_5\{(\tilde{e}_{12}+\tilde{e}_{34})\wedge \tilde{e}_5\}.
\end{align*}

\begin{thm}\label{thm_structure_deformation_SU2}
Let $(M^6,g,J)$ be an almost hermitian structure with parallel torsion related to $SU(2)$. Then there exists a $S^2$-family of orthogonal almost complex structures $\tilde{J}$ such that $(M^6,g,\tilde{J})$ is an almost hermitian structure with parallel torsion related to $SU(2)$, the corresponding characteristic connections coincide, and any orbit type of the torsion form is realized. In particular, there exists a structure of strict type $W_4$ (generalized Hopf structure).
\end{thm}

The existence of a generalized Hopf structure is of importance. 

\begin{cor}\label{cor_metric_closed_torsion_SU2}
The metric is never Einstein. 
\end{cor}

\begin{proof}
The metric and the torsion form remain unchanged when we replace the almost complex structure and consider the new almost hermitian structure. The results now follows from the fact that generalized Hopf structures are never Einstein, S. Kashiwada \cite{Kashiwada:1979}.
\end{proof}

The geometry of generalized Hopf structures is related to the geometry of Sasaki structures. A Sasaki structure on a Riemannian manifold $(X^n,g)$ of odd dimension is a Killing vector field $\xi$, the Reeb vector field, such that the endomorphism $\phi$ defined by $\phi(X)=-\nabla^g_X\xi$ satisfies $(\nabla^g_X\phi)(Y)=g(X,Y)\xi-g(\xi,Y)X$ for all vector fields $X,Y$. The 1-form dual to $\xi$ we denote by $\eta$. The characteristic connection $\nabla^c$ of a Sasaki structure $(X^n,g,\phi,\eta,\xi)$ is the unique metric connection $\nabla^c$ with totally skew-symmetric torsion that preserves the Sasaki structure, i.e. $\nabla^cg=0$, $\nabla^c\xi=0$, $\nabla^c\phi=0$, $\nabla^c\eta=0$. The torsion form, $T^c=d\eta\wedge\eta$, is parallel, $\nabla^cT^c=0$.

\begin{thm}\label{thm_generalized_hopf_SU2}
Let $(M^6,g,J)$ be a compact regular generalized Hopf structure with $Hol(\nabla^c)\subset SU(2)$. Then $M^6$ is the total space of a flat principal $S^1$-bundle over a compact five dimensional manifold $X^5$ equipped with a Sasaki structure with $Hol(\nabla^c)\subset SU(2)$. Conversely, let $M^6$ be the total space of a flat principal $S^1$-bundle over a compact five dimensional manifold $X^5$ equipped with a Sasaki structure with $Hol(\nabla^c)\subset SU(2)$. Then $M^6$ carries a generalized Hopf structure with $Hol(\nabla^c)\subset SU(2)$.
\end{thm}

\noindent
The theorem is a specialization of a general structure theorem for generalized Hopf structures by I. Vaisman, \cite[Theorem 4.1]{Vaisman:1979}. (A generalized Hopf structure is regular if the distribution determined by the vector field $JX$ is regular.)

A Sasaki structure $(X^n,g,\phi,\eta,\xi)$ is $\eta$-Einstein, if there exist functions $\lambda,\mu$ such that the Riemannian Ricci tensor is given by $Ric^g=\lambda\;g+\mu\;\eta\ox\eta$. 

\begin{prop}\label{prop_Sasaki_curvature_SU2}
A Sasaki structure $(X^5,g,\phi,\eta,\xi)$ with $Hol(\nabla^c)\subset SU(2)$ is $\eta$-Einstein. It holds
\begin{gather*}
Ric^g=6\;g-2\;\eta\ox\eta,\quad Ric^c=4\;g-4\;\eta\ox\eta.
\end{gather*}
\end{prop}

\begin{exa}\label{exa_generalized_Hopf_S5xS1}
The Hopf fibration $S^1\ra S^5\ra\CP^2$ gives rise to a Sasaki structure on $S^5$. The metric is obtained as the direct sum of the pull back metric from $\CP^2$, and the standard metric on fibers, both properly rescaled. For the Reeb vector field, one of the two unit length vector fields on the fibers may be chosen. The structure is $SU(3)$-invariant. The metric is not the round metric, but it is naturally reductive. The characteristic connection $\nabla^c$ and canonical connection $\nabla^{can}$ of $S^5=SU(3)/SU(2)$ coincide, and it follows $Hol(\nabla^c)=SU(2)$. The second cohomology with integer coefficients $H^2(S^5,\Z)=0$, and the only $S^1$-bundle over $S^5$ is the trivial bundle, $M^6=S^5\times S^1$. We obtain a generalized Hopf structure $(M^6,g,J)$ with $Hol(\nabla^c)=SU(2)$. It is not difficult to show, that any locally homogenous generalized Hopf structure with $Hol(\nabla^c)=SU(2)$ is (locally) isomorphic to it.
\end{exa}

\begin{exa}\label{exa_relatedSU2_S5xS1}
The complex structure on $S^5\times S^1$ may be replace by an almost complex structure to obtain an almost hermitian structure with parallel torsion related to $SU(2)$ realizing any orbit type of the torsion form. 
\end{exa}

\begin{NB}
If $(X^5,g,\phi,\eta,\xi)$ is an $\eta$-Einstein Sasaki structure, then, for $c>0$,
\begin{gather*}
g_c=c\;g+c(c-1)\;\eta\ox\eta,\quad\xi_c=\frac{1}{c}\xi
\end{gather*}
defines a family $(X^5,g_c,\phi_c,\eta_c,\xi_c)$ of $\eta$-Einstein Sasaki structures, see S. Tanno in \cite{Tanno:1965}. For one parameter the structure is Einstein. In a recent series of papers J. P. Gauntlett et al. presented non-homogenous Einstein Sasaki structures, see e.g. \cite{GauntlettMartelliSparksWaldram:2004}. These may be deformed to an $\eta$-Einstein Sasaki structure that satisfies the curvature conditions of Proposition \ref{prop_Sasaki_curvature_SU2}. It remains to be seen whether $Hol(\nabla^c)\subset SU(2)$. If so, these examples would yield non-homogenous almost hermitian structures with parallel torsion related to $SU(2)$.
\end{NB}

\begin{prop}\label{prop_trivialholSU2impossible}
$Hol_o(\nabla^c)=\{e\}$ is not possible.
\end{prop}

\begin{proof}
Trivial holonomy implies that $(M^6,g,J)$ locally is isomorphic to some Lie group. However, the almost complex structure may be replaced by a new almost complex structure such that the strict type is $W_4$, and Lie groups of that type do not exist, (Proposition \ref{prop_TcscalarinCliff} and the remarks before). 
\end{proof}

\noindent
It follows that $Hol_o(\nabla^c)=SU(2)$ or $T^1$. The number of parallel spinors, and eigenvalues of the action of the torsion form on the bundle of parallel spinors are given in Table \ref{tbl_par_spinors_SU2}.

\begin{NB}\label{rmk_generalized_generalized_Hopf_structures}
If the strict type is $W_1\op W_3$, then $Hol(\nabla^c)\not=Hol_o(\nabla^c)\subset SU(2)$ is possible. This class gives a slight generalization of generalized Hopf structures with $Hol(\nabla^c)\subset SU(2)$. Locally, the geometry is the same. A parallel vector field, however, need not exist globally, a parallel 1-dimensional distribution is sufficient.
\end{NB}

\begin{exa}\label{exa_generalized_generalized_Hopf_structures}
We continue with the foregoing examples. Let $\Gamma$ be the group of order 4 generated by the element $a=diag(i,i,-1)\in SU(3)$. The quotient
\begin{gather*}
M^6=S^5/\Gamma\x S^1=SU(3)/(\Gamma\x SU(2))\x S^1
\end{gather*}
is a compact manifold. The isotropy representation $\Gamma\x SU(2)\ra Aut(T_eM^6)$ at the ``origin'' of $M^6$ is the restriction of the standard representation of $SU(3)$. Now
\begin{gather*}
a\cdot(e_{12}+e_{34})\wedge e_5=-(e_{12}+e_{34})\wedge e_5
\end{gather*}
and, consequently, the generalized Hopf structure on $S^5\x S^1$ does not descend. In contrast
\begin{gather*}
a\cdot(e_{14}+e_{23})\wedge e_5=(e_{14}+e_{23})\wedge e_5.
\end{gather*}
The almost hermitian structure of strict type $W_1\op W_3$ therefore does descend, and $(M^6,g,J)$ is an almost hermitian structure with parallel torsion of strict type $W_1\op W_3$ and $Hol(\nabla^c)=\Gamma\x SU(2)$. 
\end{exa}

\begin{table}
\bdm
\begin{array}{c|c|c}
\hline
\text{holonomy}  & \text{par. spinors} & \text{eigenvalues}
\spacedhline
\hline
SU(2) & 4 & \pm\sqrt{2}\|T^c\| \\
T^1 & 4 & \pm\sqrt{2}\|T^c\|
\spacedhline
\end{array}
\edm
\caption{}
\label{tbl_par_spinors_SU2}
\end{table}

%\newpage
%----------------------------------------------------------------------------
\subsection{Structures related to $SO(3)$}\label{sec_geometry_SO3}
%----------------------------------------------------------------------------
Let $(M^6,g,J)$ be an almost hermitian structure with parallel torsion such that
\begin{gather*}
Hol(\nabla^c)\subset SO(3),\quad SO(3)\subset Iso(T^c). 
\end{gather*}
The strict type then is (i) $W_1$, (ii) $W_1\op W_3$, or (iii) $W_3$. There exists a local from such that the torsion form
\begin{align*}
T^c=\;&\alpha_1\{(e_{14}+e_{23})\wedge e_5+(e_{13}-e_{24})\wedge e_6\}\\
&+\alpha_2\{(-e_{13}+e_{24})\wedge e_5+(e_{14}+e_{23})\wedge e_6\} \\
&+\alpha_3\{3e_{135}+e_{146}+e_{236}+e_{245}\},
\end{align*}
with $\alpha_1,\alpha_2,\alpha_3$ real numbers satisfying one of the following conditions: (i) $\alpha_1>0$, $\alpha_2=\alpha_3=0$, (ii) $(\alpha_1,\alpha_2)\neq(0,0)$, $\alpha_3>0$, (iii) $\alpha_1=\alpha_2=0$, $\alpha_3>0$. Then
\begin{gather*}
\|T_2\|^2=4(\alpha_1^2+\alpha_2^2),\quad\|T_{12}\|^2=12\alpha_3^2.
\end{gather*}
(Cf. Theorem \ref{thm_strict_types_and_isotropy_groups}, and cases VII, XI, X of Table \ref{tbl_torsion_form_parameters}.)
\begin{prop}\label{lem_Ambrose_Singer_SO3}
$\nabla^cR^c=0$.
\end{prop}

\begin{proof}
The representation of $SO(3)$ on $\R^6$ is the direct sum of two real irreducible faithful representations. Hence the Berger algebra of $\so(3)$ is trivial, \cite[Lemma 4.4]{CleytonSwann:2002}, and the curvature is parallel, \cite[Lemma 5.6]{CleytonSwann:2002}. 
\end{proof}

\noindent
In fact, the almost hermitian structure $(M^6,g,J)$ is locally isomorphic to an isotropy irreducible homogeneous space $G/H$, \cite[Proposition 5.12]{CleytonSwann:2002}.

\begin{prop}\label{lem_SO3_Hol_T1_impossible}
$Hol_o(\nabla^c)=T^1$ is not possible.
\end{prop}

\begin{proof}
Computation. 
\end{proof}

Let $R^c:\Lambda^2(\R^6)\ra \so(3)$ be the curvature tensor at some point.

\begin{lem}\label{lem_curvature_so3}
$R^c=-\lambda\;pr_{\so(3)},\quad\lambda=\| T_2\|^2-\frac{1}{3}\| T_{12}\|^2$.
\end{lem}
\begin{proof}
The holonomy is trivial if and only if $\| T_2\|^2-\frac{1}{3}\| T_{12}\|^2=0$ (cf. Proposition \ref{prop_TcscalarinCliff} and the remarks before). The statement is true in that case.

The case $Hol_o(\nabla^c)=SO(3)$ remains. The curvature is parallel, and $R^c$ is $SO(3)$-equivariant. The $SO(3)$-module $\Lambda^2(\R^6)$ decomposes into three times the standard, a trivial, and a 5-dimensional module, $\Lambda^2(\R^6)=3\cdot\so(3)\op\R\op Sym^2_0(\R^3)$. A linear map $\so(3)\ra\so(3)$, or equivalently $\R^3\ra\R^3$, is $SO(3)$-equivariant if and only if it is a (real) multiple of the identity. Thus $R^c$ depends on at most three real parameters. The first Bianchi identities show that $R^c$ is completely determined by $T^c$ as stated.
\end{proof}

\noindent
The curvature tensor is determined by the torsion form, and is the same as that of $S^3\x S^3$ with its standard nearly K\"ahler structure.

A method for the classification of the local structure is outlined in Section \ref{sec_naturally_reductive_spaces}. It is readily verified, that the pairs $(R^c,T^c)$, as now determined, are infinitesimal models for all parameters. They give a complete description of the local structure. In particular, the local structure of structures of strict type $W_1$ and strict type $W_3$ is unique.

Explicit families of structures that generalize the unique structure of strict type $W_1$ on $S^3\x S^3$, and the unique structure of strict type $W_3$ on $SL(2,\C)$, are given below in Section \ref{sec_examples_SO3}. We also construct structures on $\tilde{E}_3$ and $N^6$, where $\tilde{E}_3=SU(2)\lsemidirect\R^3$ is the universal covering group of the group of euclidian motions, and $N^6$ is some nilpotent Lie group. 

It seems likely that any structure related to $SO(3)$ is locally isomorphic to a structure on one of these manifolds. The identification of the invariants $\alpha_1$, $\alpha_2$, and $\alpha_3$ in a given example, however, is very difficult. If we knew these invariants we could establish a correspondence between infinitesimal models and actual examples. Thereby we would obtain a good classification.

\begin{prop}
An almost hermitian structure related to $SO(3)$ is Einstein if and only if it is locally isometric to $S^3\x S^3$ equipped with the nearly K\"ahler or the standard product metric.
\end{prop}

\begin{proof}
A direct calculation shows, the structure is Einstein if and only if either $\alpha_1>0$, $\alpha_2=\alpha_3=0$ or $\alpha_1=0$, $\alpha_2=\alpha_3>0$. In both cases the structure is unique and is realized on $S^3\x S^3$ for the stated metrics. 
\end{proof}

\begin{prop}\label{prop_dTc_SO3}
$dT^c=\lambda\hodge\Omega,\quad\lambda=\| T_2\|^2-\frac{1}{3}\| T_{12}\|^2$.
\end{prop}

The number of parallel spinors, and eigenvalues of the action of the torsion form on the bundle of parallel spinors are given in Table \ref{tbl_par_spinors_SO3}.

\begin{table}
\bdm
\begin{array}{c|c|c}
\hline
\text{holonomy}  & \text{par. spinors} & \text{eigenvalues}
\spacedhline
\hline
SO(3)& 2 & \pm 2\|T_2\| \\
\{e\}& 8 & \pm 2\|T_2\|,\;\pm\frac{2}{\sqrt{3}}\|T_{12}\|
\spacedhline
\end{array}
\edm
\caption{}
\label{tbl_par_spinors_SO3}
\end{table}

%\newpage
%----------------------------------------------------------------------------
\subsection{Structures related to $T^2$}\label{sec_geometry_T2}
%----------------------------------------------------------------------------
\noindent
Let $(M^6,g,J)$ be an almost hermitian structure with parallel torsion such that
\begin{gather*}
Hol(\nabla^c)\subset T^2,\quad T^2\subset Iso(T^c). 
\end{gather*}
The strict type then is (i) $W_3$, (ii) $W_3\op W_4$, or (iii) $W_4$. There exists a local frame such that the torsion form
\begin{align*}
T^c=(e_{12}-e_{34})\wedge (\alpha_3e_5+\alpha_4e_6)+\alpha_5(e_{12}+e_{34})\wedge e_5
\end{align*}
with $\alpha_3,\alpha_4,\alpha_5$ real numbers satisfying one of the following conditions: (i) $\alpha_3>0$, $\alpha_4=\alpha_5=0$, (ii) $(\alpha_3,\alpha_4)\neq(0,0)$, $\alpha_5>0$, (iii) $\alpha_3=\alpha_4=0$, $\alpha_5>0$. Then
\begin{gather*}
\|T_{12}\|^2=2(\alpha_3^2+\alpha_4^2),\quad\|T_6\|^2=2\alpha_5^2.
\end{gather*}
(Cf. Theorem \ref{thm_strict_types_and_isotropy_groups}, and cases IX, IV, I of Table \ref{tbl_torsion_form_parameters}.) Throughout, let $e=(e_1,...,e_6)$ be a local frame field derived from the normal form of the torsion form. 

The tangent bundle splits
\begin{gather*}
TM^6=\mcH\op\mcV,\quad \mcH=\mcH_1\op\mcH_2,
\end{gather*}
into two horizontal distributions $\mcH_1=\{V\in TM^6:e_{12}\wedge V=0\}$, and $\mcH_2=\{V\in TM^6:e_{34}\wedge V=0\}$, and a vertical distribution $\mcV=\{V\in TM^6:e_{56}\wedge V=0\}$. The distributions $\mcH_1$,$\mcH_2$ and $\mcV$ are complex 1-dimensional and parallel.
 
\begin{prop}
%----------
The distribution $\mcV$ is integrable. Its leaves are flat and totally geodesic submanifolds.
\end{prop}

\begin{proof}
The distribution $\mcV$ is spanned by the parallel vector fields $e_5$ and $e_6$. Now $\nabla^g_{e_5}e_6=\nabla^c_{e_5}e_6-\frac{1}{2}T^c(e_5,e_6)=0$, and analogously $\nabla^g_{e_6}e_5=0$. Thus $[e_5,e_6]=\nabla^g_{e_5}e_6-\nabla^g_{e_6}e_5=0$. The proposition follows.
\end{proof}

\begin{NB}
The (local) vector fields $e_5$, $e_6$ are Killing vector fields since they are parallel with respect to a metric connection with totally skew-symmetric torsion.
\end{NB}

\noindent
Let $\Gamma(\mcH)$ and $\Gamma(\mcV)$ denote the space of (local) sections of $\mcH$ and $\mcV$.

\begin{lem}\label{lem_computation_T2}
%----------
It holds
\begin{enumerate}
\item $de_{12}=de_{34}=0$, 
\item $\Lie{V}e_{12}=\Lie{V}e_{34}=0,\quad V\in\Gamma(\mcV)$,
\item $(\nabla^g_{V_0}e_{12})(V_1,V_2)=(\nabla^g_{V_0}e_{34})(V_1,V_2)=0,\quad V_0,V_1,V_2\in\Gamma(\mcH)$.
\end{enumerate}
\end{lem}

\begin{proof}
The forms $e_{12}$ and $e_{34}$ are $\nabla^c$-parallel and (1) follows by direct calculation. The Lie derivative with respect to a vector field $V$ of a general $k$-form $\omega$ can be expressed through  $\Lie{V}\omega=V\haken d\omega+d(V\haken\omega)$ and (2) follows immediately from (1). For the covariant derivative  of a general 2-form $\omega$ we have
\begin{gather*}
(\nabla^g_{V_0}\omega)(V_1,V_2)=(\nabla^c_{V_0}\omega)(V_1,V_2)-\frac{1}{2}\omega(T^c(V_0,V_1),V_2)-\frac{1}{2}\omega(V_1,T^c(V_0,V_2))
\end{gather*}
where $V_0,V_1,V_2\in TM^6$. If $\omega$ is $\nabla^c$-parallel the first term vanishes. Furthermore $V_i,V_j\in\mcH$  implies $T^c(V_i,V_j)\in\mcV$  and (3) follows.
\end{proof}

The manifold fibers over a 4-dimensional space, 
\begin{gather*}
M^6\ra X^4=M^6/\mcV.
\end{gather*}
If $\mcV$ is a regular distribution, the quotient $X^4$ is smooth, the projection $M^6\ra X^4$ is a Riemannian submersion, and the vector fields $e_5$, $e_6$ are tangent to the vertical fibers. The case, when $M^6$ is compact and is the total space of a principle $T^2$-bundle over $X^4$, is of particular interest. In this situation, the distribution $\mcV$ defines an infinitesimal connection in $M^6$. Its curvature is determined by the differentials
\begin{gather*}
de_5=\alpha_3(e_{12}-e_{34})+\alpha_5(e_{12}+e_{34}),\quad de_6=\alpha_4(e_{12}-e_{34}).
\end{gather*}
These determine elements in the second cohomology of the base, $de_5,de_6\in H^2(X^4,\Z)$. The forms $e_{12}+e_{34}$ and $e_{12}-e_{34}$ are globally defined. They  project to $X^4$ and are parallel with respect to its Levi-Civita connection by Lemma \ref{lem_computation_T2}. This characterization of the almost hermitian structure is sufficient to reverse the construction.

\begin{thm}\label{thm_structure_T2}
%----------
Let $(M^6,g,J)$ be an almost hermitian structure with parallel torsion related to $T^2$, and assume further that $M^6$ is compact and is the total space of a $T^2$-bundle over a compact manifold $X^4$ as described. Then $X^4$ carries two K\"ahler structures of opposite orientation, and, in particular, locally is the product of two K\"ahler manifolds. Conversely, let $X^4$ be a Riemannian manifold equipped with two K\"ahler structures of opposite orientation. Then any $T^2$-bundle over $X^4$ with an infinitesimal connection as described carries a almost hermitian structure with parallel torsion related to $T^2$.
\end{thm}

\noindent
The theorem extends \cite[Theorem 4.3]{AlexandrovFriedrichSchoemann:2005}, where it was proved for structures of strict type $W_3\op W_4$.  

\begin{exa}\label{exa_T2_bundle}
Compact examples may be constructed by choosing $X^4$ to be $S^2\x S^2$, $S^2\x T^2$, $T^2\x T^2$ or, in fact, the product of any two compact K\"ahler surfaces. This also yields examples not locally homogenous.
\end{exa}

\begin{exa}
The total space of a torus bundle over a torus is a compact 2-step nilmanifold and vice versa. For the hermitian structures that arise on various nilmanifolds by choosing as base manifold the 4-dimensional torus $X^4=T^2\x T^2$ see Section \ref{sec_examples_nilmanifolds}.
\end{exa}

Let $R^c:\Lambda^2(\R^6)\ra \so(3)$ be the curvature tensor at some point.

\begin{lem}\label{lem_curvature_T2}
Let $\lambda=\frac{1}{2}(\|T_{12}\|^2-\|T_6\|^2)$. Then 
\begin{gather*}
R^c=\lambda_1\;e_{12}\ox e_{12}+\lambda_2\;e_{34}\ox e_{34}-\lambda\;(e_{12}\ox e_{34}+e_{34}\ox e_{12})
\end{gather*}
for two local functions $\lambda_1,\lambda_2$. 
\end{lem}

\begin{proof}
Follows from the Bianchi identities.
\end{proof}

\noindent
If $Hol(\nabla^c)\subset T^1$ then
\begin{gather*}
R^c=\lambda\;(e_{12}-e_{34})\ox (e_{12}-e_{34})=2\lambda\;pr_{\lan e_{12}-e_{34}\ran}.
\end{gather*}
Here $pr_{\lan e_{12}-e_{34}\ran}:\u(3)\ra\u(3)$ denotes the projection onto the Lie subalgebra generated by $e_{12}-e_{34}\in \u(3)$. (Note that $T^1=T^2\cap SU(3)$.)

\begin{exa}\label{exa_discrete_holonomy_S3_S3}
Consider the Lie group $S^3\x S^3$ with Lie algebra $\su(2)\op \su(2)$. For $i=1,2$, let $B_i(X,Y)=-\frac{1}{2}tr(XY)$, $X,Y\in\su(2)$, be the Killing form on the first and second factor respectively. For $s,t>0$, let $B_{st}=\frac{1}{s^2}B_1+\frac{1}{t^2}B_2 $. This induces a 2-parameter family of bi-invariant metrics on the manifold, which we denote by $B_{st}$ also. The Hopf fibration $S^1\ra S^3\ra S^2$ determines a distinguished direction on each sphere, and, thereby, a 2-dimensional vertical distribution $\mcV$ on $S^3\x S^3$. For $i=1,2$, let $\mcH_i$, be the intersection of the orthogonal complement of $\mcV$ with the tangent bundle of the first and second factor respectively, $T(S^3\x S^3)=\mcH_1\op\mcH_2\op\mcV$. The pull back of the complex structure of $S^2$ determines a complex structure in $\mcH_1$ and $\mcH_2$. In $\mcV$ a complex structure is defined by ``exchanging tangent vectors'', $J(X,Y)=(-Y,X)$ for $(X,Y)\in\mcV\subset T(S^3\x S^3)$. This almost complex structure is $B_{st}$-orthogonal. If $(X_0,X_1,X_2)$ is an orthonormal basis of $\su(2)$ such that $[X_i,X_{i+1}]=X_{i+2}$, $i\in\Z_3$, a right-invariant adapted frame field is given by 
\begin{gather*}
e=(e_1,...,e_6)=(s(X_0,0),s(X_1,0),t(0,X_0),t(0,X_1),s(X_2,0),t(0,X_2)). 
\end{gather*}
The torsion form of the characteristic (or canonical) connection is given by  $T^c(X,Y,Z)=-B_{st}([X,Y],Z)$ for $X,Y,Z\in T_o(S^3\x S^3)$. One finds $T^c=-s\,e_{125}-t\;e_{346}$, and thus $\|T_2\|^2=0$, $\|T_{12}\|^2=\|T_6\|^2=\frac{1}{2}(s^2+t^2)$. More precisely
\begin{gather*}
\alpha_3=\frac{s^2-t^2}{2\sqrt{s^2+t^2}},\quad\alpha_4=-\frac{st}{\sqrt{s^2+t^2}}\neq0,\quad\alpha_5=\frac{1}{2}\sqrt{s^2+t^2}.
\end{gather*}
Changing the orientation of the complex structure in the vertical distribution changes the sign of $\alpha_4$. The standard product metric is obtained for $s=t$. 
\end{exa}

\begin{exa}\label{exa_discrete_holonomy_S3_T3}
On the Lie group $S^3\x T^3$  one defines two 1-parameter families of almost hermitian structures in a similar way. One finds $\alpha_3=\pm\alpha_5$, $\alpha_4=0$.
\end{exa}

\begin{thm}\label{thm_lie_groups_T2}
An almost hermitian structure related to $T^2$ with trivial holonomy is locally isomorphic to the Lie group $S^3\x S^3$ if $\alpha_4\neq 0$, and to $S^3\x T^3$ if $\alpha_4=0$, and both Lie groups are equipped with an almost hermitian structure as described in the foregoing examples.
\end{thm}

\begin{proof}
Analogous to the proof of Theorem  \ref{thm_classification_holonomy_T1_T2} below. 
\end{proof}

\begin{exa}
We generalize the structure on $S^3\x S^3$ of Example \ref{exa_discrete_holonomy_S3_S3}. The Hopf fibration $S^1\ra S^3\ra S^2$ makes $S^3\x S^3$ into a $T^2$-principle bundle over $S^2\x S^2$. There is a natural action of $T^2\x S^3\x S^3$ on $S^3\x S^3$,
\begin{gather*}
(t_1,t_2,a_1,a_2)(b_1,b_2)=(a_1b_1t_1^{-1},a_2b_2t_2^{-1}).
\end{gather*}
Let $G=T^1\x S^3\x S^3$ with $T^1\subset T^2$ embedded anti-diagonally, i.e. $T^1=\{(t,t^{-1})\}\subset\{(t_1,t_2)\}=T^2$. Then $G$ acts on $S^3\x S^3$, and the action is effective and transitive. The isotropy group at the origin is $H=\{(t,t^{-1},t,t^{-1})\}$. Let $\alpha_3,\alpha_4,\alpha_5$ be real parameters with 
\begin{gather*}
\alpha_3+\alpha_5>0,\quad\alpha_3-\alpha_5<0,\quad\alpha_5>0.
\end{gather*}
Set $\lambda=\alpha_3^2+\alpha_4^2-\alpha_5^2$, $s=\sqrt{2\alpha_5(\alpha_3+\alpha_5)}$, $t=\sqrt{-2\alpha_5(\alpha_3-\alpha_5)}$. Let $(X_0,X_1,X_2)$ be a basis of $\su(2)$ such that $[X_i,X_{i+1}]=X_{i+2}$, $i\in\Z_3$ (as in Example \ref{exa_discrete_holonomy_S3_S3}), and such that $X_2$ is tangent to the $S^1$-fibre. Let $\g$ and $\h$ denote the Lie algebras of $G$ and $H$. Then
\begin{gather*}
\g=\R^1\op\su(2)\op\su(2)\subset\R^1\op\R^1\op\su(2)\op\su(2),\quad\h=Lin_{\R}\{(X_2,-X_2,X_2,-X_2)\},
\end{gather*}
where $\R^1\op\R^1$ is identified with a copy of the Lie algebra of the maximal torus of $S^3\x S^3$. We decompose $\g$ and define a frame $e=(e_1,...,e_6)$ by letting
\begin{gather*}
\g=\h\op\m,\quad \m=\m_+\op\m_-\op\m_0,\\
\m_+=Lin_{\R}\{e_1,e_2\},\quad\m_-=Lin_{\R}\{e_3,e_4\},\quad\m_0=Lin_{\R}\{e_5,e_6\},
\intertext{with}
e_1=(0,0,s\,X_0,0),\quad e_2=(0,0,s\,X_1,0),\quad e_3=(0,0,0,t\,X_0), \\
e_4=(0,0,0,t\,X_1),\quad e_5=-\frac{1}{2\alpha_5}(0,0,s^2\, X_2,t^2\,X_2),\\
 e_6=\frac{1}{2\alpha_4}(-2\lambda\,X_2,2\lambda\,X_2,((\frac{\alpha_3}{\alpha_5}-1)s^2-2\lambda)\,X_2,((\frac{\alpha_3}{\alpha_5}+1)t^2+2\lambda)\,X_2).
\end{gather*}
Then
\begin{align*}
[\h,\m_+]&\subset\m_+,&[\h,\m_-]&\subset\m_-,&[\h,\m_0]&=0,\\
[\m_+,\m_+]&\subset\h\op\m_0,&[\m_-,\m_-]&\subset\h\op\m_0,&[\m_+,\m_-]&=0,\\
[\m_0,\m_+]&\subset\m_+,&[\m_0,\m_-]&\subset\m_-,&[\m_0,\m_0]&=0.
\end{align*}
The commutator relations imply the existence of an invariant almost hermitian structure at once. With respect to $e,$ elements of $\h$ act on $\m$ by multiples of $e_{12}-e_{34}$. Therefore $e$ determines a family of almost hermitian structures. Geometrically, it may be characterized as the sum of the pull back of the two K\"ahler structures from the base manifold $S^2\x S^2$ and a certain metric and almost complex structure on the $T^2$-fibers. A computation shows that the metric is naturally reductive; hence $(S^3\x S^3,g,J)$ is an almost hermitian structure with parallel torsion. The torsion form of the characteristic (or canonical) connection is given by
\begin{gather*}
T^c=(e_{12}-e_{34})\wedge (\alpha_3e_5+\alpha_4e_6)+\alpha_5(e_{12}+e_{34})\wedge e_5
\end{gather*}
and the curvature by
\begin{gather*}
R^c=\lambda\;(e_{12}-e_{34})\ox (e_{12}-e_{34}).
\end{gather*}
For $\lambda=0$ the structures of Example \ref{exa_discrete_holonomy_S3_S3} are recovered.
\end{exa}

\begin{exa}
On the Lie groups $S^3\x SL(2,\R)$, $S^3\x N(1,1)$ and $T^3\x N(1,1)$ one defines families of almost hermitian structures in a similar way. Here $N(1,1)$ is the (solvable) group of rigid motions of Minkowski 2-space, i.e. the semi-direct product $\R^1\lsemidirect\R^2$ with $t\in\R^1$ acting on $\R^2$ through the matrix $\bigl[\begin{smallmatrix} e^{t} & 0 \\ 0 & e^{-t} \end{smallmatrix}\bigr]$. 
\end{exa}

\begin{thm}\label{thm_classification_holonomy_T1_T2}
An almost hermitian structure related to $T^2$ with $Hol(\nabla^c)=T^1$ is locally isomorphic to a Lie group $G=(T^1\x G)/T^1$ as given in Table \ref{tbl_local_models_T2_T1}, equipped with an invariant structure as described in the foregoing examples.
\end{thm}

\begin{proof}
The method of proof is outlined in Section \ref{sec_naturally_reductive_spaces}. (The pairs $(R^c,T^c)$ as determined by the classification of the torsion orbits and the remark after Lemma \ref{lem_curvature_T2} is an infinitesimal model for all parameters. One verifies that all infinitesimal models are realized on the given manifolds as stated.)
\end{proof}

\begin{table}
\bdm
\begin{array}{c|c}
\hline
\text{orbit} & \text{Lie group}
\spacedhline
\hline
\alpha_3+\alpha_5<0 & S^3\x SL(2,\R) \\
\alpha_3+\alpha_5>0,\; \alpha_3-\alpha_5<0 & S^3\x S^3 \\
\alpha_3-\alpha_5>0 & S^3\x SL(2,\R) \\
\alpha_3=\pm\alpha_5 & S^3\x N(1,1) \\
\alpha_5=0 & T^3\x N(1,1) 
\spacedhline
\end{array}
\edm
\caption{}
\label{tbl_local_models_T2_T1}
\end{table}

\begin{NB}
A classification of the locally homogenous almost hermitian structures related to $T^2$ with holonomy $T^2$ or $T^1\neq T^2\cap SU(3)$ can be worked out by generalizing the foregoing examples. These structures, however, do not admit parallel spinors.
\end{NB}

\begin{prop}
An almost hermitian structure related to $T^2$ is Einstein if and only if it is locally isometric to $S^3\x S^3$ with the standard product metric.
\end{prop}

\begin{proof}
A direct calculation shows, the Riemannian Ricci tensor is given by the following three $2\x2$ blocks arranged along the diagonal: 
$\bigl[\begin{smallmatrix} A+\lambda_1 & 0 \\ 0 & A+\lambda_1 \end{smallmatrix}\bigr]$, 
$\bigl[\begin{smallmatrix} A+\lambda_2 & 0 \\ 0 & A+\lambda_2 \end{smallmatrix}\bigr]$ and 
$\bigl[\begin{smallmatrix} -3(\alpha_3^2+\alpha_5^2) & -3\alpha_3\alpha_4 \\ -3\alpha_3\alpha_4 & -3\alpha_4^2 \end{smallmatrix}\bigr]$, with $A=-\frac{3}{2}((\alpha_3+\alpha_5)^2+\alpha_4^2)$. Thus Einstein implies $\lambda_1=\lambda_2=0$ and $\alpha_3=0$, $\alpha_4^2=\alpha_5^2$. The curvature vanishes, $R^c=0$. Now apply Theorem \ref{thm_lie_groups_T2}.
\end{proof}

The number of parallel spinors, and eigenvalues of the action of the torsion form on the bundle of parallel spinors are given in Table \ref{tbl_par_spinors_T2}.

\begin{table}
\bdm
\begin{array}{c|c|c}
\hline
\text{holonomy}  & \text{par. spinors} & \text{eigenvalues}
\spacedhline
\hline
SU(2) & 4 & \pm\sqrt{2}\|T_6\| \\
T^1 & 4 & \pm\sqrt{2}\|T_6\|\\
\{e\} & 8 & 0,\;\pm\sqrt{2}\|T_6\|,\;\pm\sqrt{2}\|T_{12}\| 
\spacedhline
\end{array}
\edm
\caption{}
\label{tbl_par_spinors_T2}
\end{table}

%----------------------------------------------------------------------------
\subsection{Structures related to $T^1$}\label{sec_geometry_T1}
%----------------------------------------------------------------------------
Let $(M^6,g,J)$ be an almost hermitian structure with parallel torsion such that strict type, holonomy and isotropy group satisfy one of following conditions: (i) $W_1\op W_3$, $Iso_o(T^c)=T^1$,  (ii) $W_1\op W_3\op W_4$, $Iso(T^c)=T^1$, (iii) $W_3\op W_4$, $Hol(\nabla^c)\subset T^1$.  There exists a local frame such that the torsion 
\begin{gather*}
T^c=\alpha_1\{(e_{14}+e_{23})\wedge e_5\}+(e_{12}-e_{34})\wedge (\alpha_3e_5+\alpha_4e_6)+\alpha_5\{(e_{12}+e_{34})\wedge e_5\},
\end{gather*}
with $\alpha_1,\alpha_3,\alpha_4,\alpha_5$ real numbers satisfying on of the following conditions:  (i)  $\alpha_1>0$, $\alpha_5=0$, (ii)  $\alpha_1>0$, $\alpha_5>0$, (iii)  $\alpha_1=0$, $\alpha_5>0$ , and in all cases either $\alpha_3>0$, $\alpha_4\in\R$ or $\alpha_3=0$, $\alpha_4>0$. Then
\begin{gather*}
\|T_2\|^2=\alpha_1^2,\quad\|T_{12}\|^2=\alpha_1^2+2(\alpha_3^2+\alpha_4^2),\quad\|T_6\|^2=2\alpha_5^2.
\end{gather*}
(Cf. Theorem \ref{thm_strict_types_and_isotropy_groups}, and cases III, VI, IV of Table \ref{tbl_torsion_form_parameters}.) We will content ourselves with a brief analysis. 

It holds $Hol_o(\nabla^c)\subset T^1$. There exists, analogously to the case of structures related to $SU(2)$, an associated (local) $S^2$-family of parallel orthogonal almost complex structures. The corresponding almost hermitian structures, it is checked, are structures related to $T^1$. Not every orbit may be realized by replacing the almost complex structure. There exists, however, a structure of strict type $W_3\op W_4$ with $T^2\subset Iso(T^c)$. The (local) geometry of structures related to $T^1$, therefore, may be derived from the geometry of structures related to $T^2$. In particular, since the Berger algebra of $T^1$ is trivial, a structure related to $T^1$ is locally isomorphic to an isotropy irreducible homogeneous space $G/H$, and the (local) classification of structures related to $T^2$ with holonomy contained in $T^1$ yields a complete (local) classification of structures related to $T^1$ (cf. Theorem \ref{thm_lie_groups_T2} and Theorem \ref{thm_classification_holonomy_T1_T2}).
%
%---------------------------------------------------------------------------- 
\subsection{Examples of structures related to $SO(3)$}\label{sec_examples_SO3}\label{page_examples_SO3}
%----------------------------------------------------------------------------
\noindent
In the following, we define families of almost hermitian structures with parallel torsion on the manifolds $S^3\times S^3$, $SL(2,\C)$, $\tilde{E}_3$ and $N^6$, where $\tilde{E}_3=SU(2)\lsemidirect\R^3$ is the universal covering group of the group of euclidian motions of $\R^3$, and $N^6$ is some nilpotent Lie group. 

Let $M^6$ be one of the manifolds $S^3\times S^3$, $SL(2,\C)$, $\tilde{E}_3$ or $N^6$. In each case then, we choose a simply connected, connected Lie group $G$ that acts effectively and transitively on $M^6$. We determine the isotropy group $H$ at the origin of $M^6$ and thus obtain a homogenous structure, $M^6\isomorph G/H$. We denote by $\g$ and $\h$ the Lie algebras of $G$ and $H$. We fix a decomposition  $\g=\h\op\m$, $\m=\m_1\op\m_2$, such that, for $i=1,2$, it holds $dim(\m_i)=3$, and $[\h,\m_i]\subset\m_i$ (the decomposition is reductive). We compute the isotropy representation $\rho:H\ra Aut(\m)$, and define on $\m$ a $H$-invariant hermitian structure $(g,J)$. This defines on $M^6$ an invariant almost hermitian structure. In all cases the metric is naturally reductive, i.e. $g([X,Y]_{\m},Z)+g(Y,[X,Z]_{\m})=0$ for any $X,Y,Z\in\m$. The canonical connection is a hermitian connection with totally skew-symmetric torsion, and thus coincides with the characteristic connection. The torsion form is parallel. 

In each case the automorphism group $G$ is chosen such that the isotropy group is isomorphic to $SU(2)$, and the isotropy representation is reducible and decomposes into two real irreducible 3-dimensional representations (with kernel $\{\pm1\}$). The almost hermitian structure therefore is of type $W_1\op W_3$, and $SO(3)\subset Iso(T^c)$.

\begin{exa}\label{exa_SO3_S3xS3}
Let $G=S^3\x S^3\x S^3$ act on $M^6=S^3\x S^3$ by 
\begin{gather*}
(a_0,a_1,a_2)(b_1,b_2)=(a_1b_1a_0^{-1},a_2b_2a_0^{-1}). 
\end{gather*}
The action is effective and transitive. The isotropy group at the origin is $H=\{(h,h,h): h\in SU(2)\}$. Let $b,d$ and $k_1,k_2$ be real parameters such that $b\neq d,\quad k_1,k_2>0$. Set $a=-\frac{(d-1)(d\,k_1+b\,k_2)}{(b-d)k_2}$, $c=\frac{(b-1)(d\,k_1+b\,k_2)}{(b-d)k_1}$, and restrict to parameters $(b,d,k_1,k_2)$ such that $det\;(
\bigl[\begin{smallmatrix} 
1 & 1 & 1 \\ 
1 & a & b \\ 
1 & c & d 
\end{smallmatrix}\bigr]
)=ad+b+c-a-bc-d\neq0$. We decompose $\g=\h\op\m_1\op\m_2$, where 
\begin{align*}
\g &= \{(A,B,C): A,B,C\in\su(2)\},& \h &= \{(A,A,A): A\in\su(2)\}, \\
\m_1 &= \{(A,aA,bA): A\in\su(2)\}, & \m_2 &= \{(B,cB,dB): B\in\su(2)\}.
\end{align*}
The isotropy representation is given by 
\begin{align*}
\rho(h,h,h)&((A,aA,bA)+(B,cB,dB)) \\
&=((hAh^{-1},a\,hAh^{-1},b\,hAh^{-1})+(hBh^{-1},c\,hBh^{-1},d\,hBh^{-1})).
\end{align*}
Let $B$ denote the Killing form on $\su(2)$, $B(X,Y)=-\frac{1}{2}tr(XY)$. We define 
\begin{gather*}
g((A_1,aA_1,bA_1)+(B_1,cB_1,dB_1),(A_2,aA_2,bA_2)+(B_2,cB_2,dB_2))\\=k_1\,B(A_1,A_2)+k_2\,B(B_1,B_2), \\
J((A,aA,bA)+(B,cB,dB))=-\sqrt{\frac{k_2}{k_1}}(B,aB,bB)+\sqrt{\frac{k_1}{k_2}}(A,cA,dA).
\end{gather*}
It is not difficult to see now, that $(S^3\x S^3,g,J)$ is an almost hermitian structure with parallel torsion of type $W_1\op W_3$, and $SO(3)\subset Iso(T^c)$. Moreover,
\begin{gather*}
\|T_2\|^2=\frac{k_1+k_2}{(b-d)^2}(\frac{(d-1)^2}{k_2}+\frac{(b-1)^2}{k_1})(\frac{d^2}{k_2}+\frac{b^2}{k_1}),\quad\|T_2\|^2-\frac{1}{3}\|T_{12}\|^2=-8(\frac{b}{k_1}+\frac{d}{k_2}).
\end{gather*}
\noindent
Certain parameters are special: 
\begin{enumerate}
\item[(i)] For $b=-2$, $d=0$, $k_1=3$, $k_2=1$ it holds $T_2=\frac{16}{3}$ and $T_{12}=0$, that is, the structure is nearly K\"ahler, non-K\"ahler.
\item[(ii)] For $b=-d\frac{k_1}{k_2}$ it holds $\|T_2\|^2-\frac{1}{3}\|T_{12}\|^2=0$. The characteristic curvature vanishes and the torsion form is closed, cf. Lemma \ref{lem_curvature_so3} and Proposition \ref{prop_dTc_SO3}. 
\item[(iii)] For $b=-d=\pm1$ and $k_1=k_2$ it holds $\|T_2\|^2-\frac{1}{3}\|T_{12}\|^2=0$, and the metric is Einstein. This is the standard product metric.
\end{enumerate}
\end{exa}
 
\begin{exa}
Let $G=SU(2)\times SL(2,\C)$ act on $M^6=SL(2,\C)$ by 
\begin{gather*}
(a_0,a_1)(b_1)=(a_1 b_1 a_0^{-1}). 
\end{gather*}
The action is effective and transitive. The isotropy group at the origin is $H=\{(h,h): h\in SU(2)\}$. Let $\su(2)^{\perp}$ denote the orthogonal complement of $\su(2)$ inside $\slla(2,\C)$ with respect to the Killing form, an isomorphism $\su(2)\isomorph\su(2)^{\perp}$ being given by multiplication with $i\in\C$. We decompose $\g$ for the real parameter $p>0$,
\begin{align*}
\g &= \{(A_0,A_1+B): A_0,A_1\in\su(2),B\in\su(2)^{\perp}\}, & \h &= \{(A,A): A\in\su(2)\}, \\
\m_1 &= \{(A,\frac{1}{p+1}A): A\in\su(2)\},& \m_2 &= \{(0,B): B\in\su(2)^{\perp}\}.
\end{align*}
The isotropy representation is given by 
\begin{align*}
& \rho(h,h)(A,\frac{1}{p+1}A+B)=(hAh^{-1},\frac{1}{p+1}hAh^{-1}+hBh^{-1}). 
\end{align*}
Let $B$ denote the Killing form on $\slla(2,\C)$, $B(X,Y)=-\frac{1}{2}tr(XY)$. We define 
\begin{gather*}
g((A_1,\frac{1}{p+1}A_1+B_1),(A_2,\frac{1}{p+1}A_2+B_2))= B(A_1,A_2)-\frac{(p+1)^2}{p}B(B_1,B_2), \\
J((A,\frac{1}{p+1}A+B)) = \frac{p+1}{\sqrt{p}}(i B,\frac{p}{(p+1)^2}i A+\frac{1}{p+1}i B).
\end{gather*}
\noindent
It is not difficult to see now, that $(SL(2,\C),g,J)$ is an almost hermitian structure with parallel torsion of type $W_1\op W_3$, and $SO(3)\subset Iso(T^c)$. Moreover,
\begin{gather*}
\|T_2\|^2=\frac{(p-1)^2}{(p+1)^2},\quad\|T_{12}\|^2=\frac{3(p+3)^2}{(p+1)^2},\quad\|T_2\|^2-\frac{1}{3}\|T_{12}\|^2=-\frac{8}{p+1}.
\end{gather*}
The structure is of strict type $W_3$ if $p=1$, and of strict type $W_1\op W_3$ otherwise. (Natural reductivity may be established by using the fact that the adjoint action of an element of a semi-simple Lie algebra is skew-symmetric with respect to the Killing form.)
\end{exa}

\begin{exa}
Let $G=SU(2)\times(SU(2)\lsemidirect\R^3) $ act on $M^6=SU(2)\lsemidirect\R^3$ by
\begin{gather*}
(a_0,a_1,v)\cdot(b_1,w)=(a_1b_1a_0^{-1},a_1w+v).
\end{gather*}
The action is effective and transitive. The isotropy group at the origin is $H=\{(h,h,0): h\in SU(2)\}$. We decompose $\g$, 
\begin{align*}
\g &= \{(A_0,A_1,v): A_0,A_1\in\su(2), v\in\R^3\}, & \h &= \{(A,A,0): A\in\su(2)\}, \\
\m_1 &= \{(0,0,v): A\in\su(2), v\in\R^3\}, & \m_2 &= \{(A,0,0): A\in\su(2), v\in\R^3\}.
\end{align*}
The isotropy representation is given by 
\begin{align*}
\rho(h,h,0)(A,0,v)=(hAh^{-1},0,hv).
\end{align*}
Let $\lan,\ran$ denote the standard positive definite inner product on $\R^3$ and $\hodge$ the Hodge operator on $(\R^3,\lan,\ran)$. We identify $\su(2)=\so(3)$ and define  
\begin{gather*}
g((A_1,0,v_1),(A_2,0,v_2)) = \lan v_1,v_2\ran+\lan\hodge A_1,\hodge A_2\ran, \\
J((A,0,v)) = (-\hodge v,0,\hodge A). 
\end{gather*}
\noindent
It is not difficult to see now, that $(\tilde{E}_3,g,J)$ is an almost hermitian structure with parallel torsion of type $W_1\op W_3$, and $SO(3)\subset Iso(T^c)$. Moreover,
\begin{gather*}
\|T_2\|^2=\frac{1}{4},\quad\|T_{12}\|^2=\frac{3}{4},\quad\|T_2\|^2-\frac{1}{3}\|T_{12}\|^2=0.
\end{gather*}
The characteristic curvature vanishes and the torsion form is closed. 
\end{exa}

\begin{exa}\label{exa_SO3_N6}
Let $N^6=\R^3\times\R^3$ be the Lie group diffeomorphic to $\R^3\times\R^3$ with multiplication $(a_1,b_1)(a_2,b_2)=(a_1+a_2,b_1+b_2+\frac{1}{2}(a_1\times a_2))$. Here $\times$ denotes the vector product of $\R^3$. The bracket of the Lie algebra $\n\isomorph\R^3\op\R^3$ is given by 
\begin{align*}
[(v_1,w_1),(v_2,w_2)]=(0,v_1\times v_2),
\end{align*}
and, consequently, $N^6$ is a 2-step nilpotent.

Furthermore, let $G=SU(2)\lsemidirect N^6$ be the Lie group diffeomorphic to $SU(2)\times N^6$ with multiplication $(h_1,a_1,b_1)(h_2,a_2,b_2)=(h_1h_2,h_1a_2+a_1,h_1b_2+b_1+\frac{1}{2}(a_1\times h_1a_2))$, the action of $SU(2)$ on $\R^3$ by the double cover $SU(2)\ra SO(3)$. (Obviously, $SU(2)\subset G$ is a subgroup, $N^6\subset G$ is normal, being the kernel of an obvious homomorphism $G\ra SU(2)$, and $SU(2)\cap N^6=\{e\}$. Therefore $G$ is a semi-direct product of $SU(2)$ and $N^6$.) The bracket of the Lie algebra $\g\isomorph\su(2)\op\n$ is given by 
\begin{align*}
[(A_1,v_1,w_1),(A_2,v_2,w_2)]=([A_1,A_2],A_1v_2-A_2v_1,A_1w_2-A_2w_1+v_1\times v_2).
\end{align*}

Let an element of $G=SU(2)\lsemidirect N^6$ act on $M^6=N^6$ by left multiplication of its projection onto $N^6$, 
\begin{gather*}
(h_1,a_1,b_1)(a_2,b_2)=(a_1,b_1)(a_2,b_2)=(a_1a_2,b_1b_2). 
\end{gather*}
The action is effective and transitive. The isotropy group at the origin is $H=\{(h,0,0):h\in SU(2)\}$. Let $\lan,\ran$ denote the standard positive definite inner product on $\R^3$, and $\hodge$ the Hodge operator on $(\R^3,\lan,\ran)$. We identify $\su(2)=\so(3)$, and decompose $\g$, 
\begin{align*}
\g&=\{(A,v,w):A\in \su(2),v,w\in\R^3\}, & \h&=\{(A,0,0):A\in \su(2)\}, \\
\m_1&=\{(0,v,0):v\in\R^3\}, & \m_2&=\{(*w,0,w):w\in\R^3\}. 
\end{align*}
The isotropy representation is given by 
\begin{align*}
\rho(h,0,0)(\hodge w,v,w)=(\hodge hw,hv,hw).
\end{align*} 
We define 
\begin{gather*}
g((*w_1,v_1,w_1),(*w_2,v_2,w_2))=\lan v_1,v_2\ran+\lan w_1,w_2\ran, \\
J((*w,v,w))=(*v,-w,v).
\end{gather*}
\noindent
It is not difficult to see now, that $(N^6,g,J)$ is an almost hermitian structure with parallel torsion of type $W_1\op W_3$, and $SO(3)\subset Iso(T^c)$. Moreover,
\begin{gather*}
\|T_2\|^2=\frac{1}{4},\quad\|T_{12}\|^2=\frac{27}{4},\quad\|T_2\|^2-\frac{1}{3}\|T_{12}\|^2=-2\neq0.
\end{gather*}
\noindent
For a treatment of this example from the point of view of invariant structures on Lie groups (rather than naturally reductive spaces) see P. Ivanov and S. Ivanov in \cite{IvanovIvanov:2004}.
\end{exa}

%----------------------------------------------------------------------------
\subsection{Examples of structures on nilmanifolds}\label{sec_examples_nilmanifolds}
%----------------------------------------------------------------------------
\noindent
A nilmanifold is a quotient $M^n=\Gamma\backslash N^n$, where $N^n$ is a simply connected nilpotent Lie group of dimension $n$, and $\Gamma$ a cocompact discrete subgroup. Nilmanifolds are compact by definition. If the structure constants of $N^n$ are rational such $\Gamma$ always exist. The structure constants of all nilpotent Lie groups of dimension $\leq6$ are rational. While in general structures with parallel torsion on nilmanifolds are rare, here, we will define families of such structures on three 2-step nilpotent Lie groups, that give rise to such structures on nilmanifolds. 

The nilpotency of a nilpotent Lie group $N^n$ is equivalent to the existence of a global frame of left-invariant 1-forms $(e_1,...,e_n)$ such that
\begin{gather*}
de_i\in\Lambda^2(e_1,...,e_{i-1}),\quad i=1,...,n.
\end{gather*}
Here, $\Lambda^2(e_1,...,e_{i-1})$ denotes the second exterior power of the vector space spanned by $(e_1,...,e_{i-1})$ for $i=2,...,n$, and the zero vector space for $i=1$. A left-invariant structure on $N^n$ descends to a structure on $M^n=\Gamma\backslash N^n$, and for the specification of an almost hermitian structure on $N^n$ the definition of a global left-invariant adapted frame is sufficient. 

Let $\alpha_3,\alpha_4,\alpha_5$ be three real parameters subject to one of the following mutually exclusive conditions: 
\begin{gather*}
(i)\;\alpha_3=\pm\alpha_5,\;\alpha_4=0,\;\alpha_5>0, \quad(ii)\;\alpha_3\neq\pm\alpha_5,\;\alpha_3\neq0,\;\alpha_4=0,\;\alpha_5>0, \\
(iii)\;\alpha_3\neq0,\;\alpha_4\neq0,\; \alpha_5>0, \quad(iv)\;\alpha_3=0,\;\alpha_4\neq0,\;\alpha_5>0,\\
(v)\;\alpha_3=\alpha_4=0,\;\alpha_5>0,\quad (vi)\;\alpha_3>0,\;\alpha_4=\alpha_5=0.
\end{gather*}
Let $N^6$ be the nilpotent Lie group with left-invariant vector fields $(e_1,...,e_6)$ whose dual 1-forms, denoted by the same symbols, satisfy the structure equations
\begin{gather*}
de_1=de_2=de_3=de_4=0,\quad de_5=\alpha_3(e_{12}-e_{34})+\alpha_5(e_{12}+e_{34}), \\
de_6=\alpha_4(e_{12}-e_{34}).
\end{gather*}
For $i=1,...,6$, obviously $dde_i=0$, and the equations define Lie groups. The nilpotency is obvious. Now, $de_i(X,Y)=-e_i([X,Y])$ for left-invariant vector fields $X,Y$, and the only non-trivial brackets are
\begin{gather*}
[e_1,e_2]=-[e_2,e_1]=-(\alpha_3+\alpha_5)e_5-\alpha_4e_6,\\
[e_3,e_4]=-[e_4,e_3]=(\alpha_3-\alpha_5)e_5-\alpha_4e_6.
\end{gather*}
Let $\n$ denote the Lie algebra of $N^6$. Then $[\n,[\n,\n]]=0$. Hence $N^6$ is 2-step nilpotent.

We compute the Nijenhuis tensor $N$. For instance
\begin{align*}
N(e_1,e_3)&=[Je_1,Je_3]-J[Je_1,e_3]-J[e_1,Je_3]-[e_1,e_3]\\
&=[e_2,e_4]-J[e_2,e_3]-J[e_1,e_4]-[e_1,e_2]=0.
\end{align*}
Analogous computations show $N(e_1,e_5)=N(e_3,e_5)=0$. The symmetries of the Nijenhuis tensor imply $N=0$. The almost hermitian structure therefore is of type $W_3\op W_4$ for all parameters, and the characteristic connection $\nabla^c$ always exists.

Let $\Omega$ denote the K\"ahler form. We compute
\begin{align*}
d\Omega&=d(e_{12}+e_{34}+e_{56})=de_5\wedge e_6-e_5\wedge de_6 \\
&=(e_{12}-e_{34})\wedge (\alpha_3e_6-\alpha_4e_5)+\alpha_5(e_{12}+e_{34})\wedge e_6.
\end{align*}
Decomposing $d\Omega$ into its $U(3)$-components, and using Theorem \ref{thm_type_and_diff_eqs}, we find
\begin{align*}
T^c=(e_{12}-e_{34})\wedge (\alpha_3e_5+\alpha_4e_6)+\alpha_5(e_{12}+e_{34})\wedge e_5.
\end{align*}
If we let $T^c=T_{12}+T_6$ as usual, then
\begin{align*}
\|T_{12}\|^2=2(\alpha_3^2+\alpha_4^2),\quad\|T_6\|^2=2\alpha_5^2.
\end{align*}
Hence the strict type is as given in Table \ref{tbl_examples_nilmanifolds}. 

If the torsion is parallel, $dT^c=\sigma_{T^c}$. This integrability condition is easy to check,
\begin{align*}
dT^c=-2(\alpha_3^2+\alpha_4^2-\alpha_5^2)e_{1234}=2\sigma_{T^c}.
\end{align*}
(Cf. page \pageref{page_sigmaTc} for the definition of $\sigma_{T^c}$.) The computation shows that the families (i), (iii) and (iv) contain structures with closed torsion, $dT^c=0$.

The verification of $\nabla^cT^c=0$ is routine, but involved, for it requires the computation of the Levi-Civita and the characteristic connection. The Levi-Civita connection is determined by
\begin{align*}
2g(\nabla^g_XY,Z)=g([XY],Z)-g([YZ],X)+g([ZX],Y).
\end{align*}
for left-invariant vector fields $X,Y,Z$. The non-trivial $\nabla^g_{e_i}e_j$ are
\begin{gather*}
\nabla^g_{e_1}e_2=-\nabla^g_{e_2}e_1=\frac{1}{2}(-(\alpha_3+\alpha_5)e_5-\alpha_4e_6), \\
\nabla^g_{e_3}e_4=-\nabla^g_{e_4}e_3=\frac{1}{2}((\alpha_3-\alpha_5)e_5-\alpha_4e_6), \\
\nabla^g_{e_1}e_5=\nabla^g_{e_5}e_1=\frac{1}{2}(\alpha_3+\alpha_5)e_2,\quad\nabla^g_{e_1}e_6=\nabla^g_{e_6}e_1=\frac{1}{2}\alpha_4e_2,\\
\nabla^g_{e_2}e_5=\nabla^g_{e_5}e_2=-\frac{1}{2}(\alpha_3+\alpha_5)e_1,\quad\nabla^g_{e_2}e_6=\nabla^g_{e_6}e_2=-\frac{1}{2}\alpha_4e_1,\\
\nabla^g_{e_3}e_5=\nabla^g_{e_5}e_3=-\frac{1}{2}(\alpha_3-\alpha_5)e_4,\quad\nabla^g_{e_3}e_6=\nabla^g_{e_6}e_3=-\frac{1}{2}\alpha_4e_4,\\
\nabla^g_{e_4}e_5=\nabla^g_{e_5}e_4=\frac{1}{2}(\alpha_3-\alpha_5)e_3,\quad\nabla^g_{e_4}e_6=\nabla^g_{e_6}e_4=\frac{1}{2}\alpha_4e_3.
\end{gather*}
For the characteristic connection
\begin{align*}
\nabla^c_{e_i}e_j=\nabla^g_{e_i}e_j+\frac{1}{2}T^c(e_i,e_j)
\end{align*}
one finds correspondingly
\begin{gather*}
\nabla^c_{e_5}e_1=(\alpha_3+\alpha_5)e_2,\quad\nabla^c_{e_6}e_1=\alpha_4e_2, \\
\nabla^c_{e_5}e_2=-(\alpha_3+\alpha_5)e_1,\quad\nabla^c_{e_6}e_2=-\alpha_4e_1, \\
\nabla^c_{e_5}e_3=(-\alpha_3+\alpha_5)e_4,\quad\nabla^c_{e_6}e_3=-\alpha_4e_4, \\
\nabla^c_{e_5}e_4=(\alpha_3-\alpha_5)e_3,\quad\nabla^c_{e_6}e_4=\alpha_4e_3.
\end{gather*}
(There are two distributions, a horizontal distribution $\mcH$ spanned by the vector fields $e_1,...,e_4$ and a vertical distribution $\mcV$ spanned by the vector fields $e_5, e_6$, and $TM^6=\mcH\op\mcV$. The vector fields $e_1,...,e_4$ are covariantly constant in any horizontal direction, their covariant derivative in a vertical direction equals a rotation (=application of $J$) times a scalar. The vector fields $e_5, e_6$ are parallel.) Now,
\begin{gather*}
(\nabla^c_XT^c)(Y,Z,U)=-T^c(\nabla^c_XY,Z,U)-T^c(Y,\nabla^c_XZ,U)-T^c(Y,Z,\nabla^c_XU)
\end{gather*}
for left-invariant vector fields $X,Y,Z,U$, and the verification of $\nabla^cT^c=0$ is straight forward.

The Betti numbers of a nilmanifold $M^6=\Gamma\backslash N^6$ equal the Betti numbers of the finite dimensional differential complex
\begin{gather*}
0\ra\n^*\xrightarrow{d}\Lambda^2(\n^*)\ra...\ra\Lambda^6(\n^*)\ra0
\end{gather*}
where $\n^*$ denotes the vector space dual to $\n$, and $d$ the exterior derivative. The Betti numbers of our examples are readily computed. 

The appendix of \cite{Salamon:2000} contains a classification of 6-dimensional nilpotent Lie algebras. There are six isomorphism classes of 2-step nilpotent Lie algebras. The pairs $b_1=5$, $b_2=11$ and $b_1=5$, $b_2=9$ each determine a unique Lie algebra. There are, however, three nilpotent Lie algebras with $b_1=4$, $b_2=8$. To identify the structure of the nilpotent Lie algebra of our examples in that case, we note, that it is not difficult to see that there exists a basis such that the structure equations are given by 
\begin{gather*}
de_1=de_2=de_3=de_4=0,\quad de_5=e_{12},\quad de_6=e_{34}, 
\end{gather*}
which in the standard short hand notation reads as $(0,0,0,0,12,34)$.

\begin{table}
\bdm
\begin{array}{c|c|c|c|c}
\hline
\text{family} & \text{strict type} & b_1 & b_2 & \text{structure}
\spacedhline
\hline
(i) & W_3\op W_4 & 5 & 11 & (0,0,0,0,0,12) \\
(ii), (iii) & W_3\op W_4 & 4 & 8 & (0,0,0,0,12,34) \\
(iv) & W_3\op W_4 & 5 & 9 & (0,0,0,0,0,12+34) \\
(v) & W_4 & 5 & 9 & (0,0,0,0,0,12+34) \\
(vi) & W_3 & 5 & 9 & (0,0,0,0,0,12+34)
\spacedhline
\end{array}
\edm
\caption{}
\label{tbl_examples_nilmanifolds}
\end{table}

%\newpage
%---------------------------------------------------------------------------- 
\subsection{Naturally reductive spaces and Lie groups}\label{sec_naturally_reductive_spaces}
%---------------------------------------------------------------------------- 
Naturally reductive spaces and Lie groups are important examples of almost hermitian structures with parallel torsion. In the following, we give an account on their classification, and put into this perspective some of the previous geometric results.

The canonical connection of a naturally reductive space is the characteristic connection of any invariant almost hermitian structure that may exist on that space. In this situation, $\nabla^c R^c=0$, $\nabla^c T^c=0$, and the classification may be approached through the classification of the infinitesimal models. If $H$ is a holonomy group, $\h$ its Lie algebra, this means the classification of pairs $(R^c,T^c)$, where $R^c:\Lambda^2(\R^n)\ra\h$ is a curvature tensor, and $T^c:\Lambda^3(\R^n)\ra\R$ is a torsion form. There are several identities any such pair satisfies, most notably, $R^c$ and $T^c$ must be $H$-invariant, and the the first Bianchi identities must hold. 

The Nomizu construction associates  Lie algebras $\g$ and $\h$ with a given infinitesimal model $(R^c,T^c)$. Briefly, let $\h=\{A\in\u(3):A.R^c=0,A.T^c=0\}$, let $\g=\h\op\R^6$, and define a bracket on $\g$ by
\begin{gather*}
[(A,X),(B,Y)]=([A,B]-R(X,Y),AX-BY-T(X,Y))
\end{gather*}
for $(A,X),(B,Y)\in\g$. There exist Lie groups $G$ and $H$ with Lie algebras $\g$ and $\h$. The space $M^n=G/H$ then is a naturally reductive space with infinitesimal model isomorphic to $(R^c,T^c)$, if and only if $H\subset G$ is closed. For details with regard to infinitesimal models and the Nomizu construction we refer to the extensive literature (e.g. \cite{Tricerri:1992}).

A preliminary discussion of the holonomy groups of the characteristic connection may be derived from a discussion of the closed, connected subgroups $H\subset Iso(T^c)$ of the various isotropy groups. By construction of examples one may show that some (or all) of these groups occur. In the following, we treat trivial, abelian, and non-abelian holonomy groups, with particular emphasis on groups that admit parallel spinors (i.e. $H\subset SU(3)$).

%\newpage
%---------------------------------------------------------------------------- 
\subsubsection{Trivial holonomy (Lie groups)}
%---------------------------------------------------------------------------- 
Let $\R^6=(\R^6,g,J)$ be the standard hermitian vector space, let $T^c$ be a torsion form. The multiplication defined on $\g=(\R^6,T^c)$ by $[X,Y]=-T^c(X,Y)$ for $X,Y\in\R^6$ does not satisfy the Jacobi identities in general. The Jacobi identities hold if and only if $(T^c)^2\in Cl(\R^6)$ is a scalar, (cf. Remark \ref{rmk_Sternberg_criteria}).

\begin{prop}\label{prop_TcscalarinCliff}
Let $(M^6,g,J)$ be an almost hermitian structure with parallel torsion. Then the following are equivalent:
\begin{enumerate}
\item[(i)] $(T^c)^2\in\R\subset Cl(\R^6)$,
\item[(ii)] $3\|T_2\|^2-\|T_{12}\|^2+\|T_6\|^2=0$,
\item[(iii)] $dT^c=0$.
\end{enumerate}
\end{prop}
\noindent

\begin{thm}\label{thm_Lie_groups}
Let $(M^6,g,J)$ be a Lie group with an invariant almost hermitian structure with parallel torsion. Then the strict type is either $W_1\op W_3$, $W_3\op W_4$ or $W_1\op W_3\op W_4$. More precisely,
\begin{enumerate}
\item[(i)] There exists an $S^1$-family of structures of strict type $W_1\op W_3$. In that case, $Iso_o(T^c)=SO(3)$. 
\item[(ii)] There exists an $S^1$-family of structures of strict type $W_3\op W_4$. In that case the underlying Lie group is isomorphic to $S^3\x S^3$ or $S^3\x T^3$.
\item[(iii)] There exists a 2-parameter family of structures on Lie groups of strict type $W_1\op W_3\op W_4$. In that case, $Iso_o(T^c)\subset T^1$, and the almost complex structure may be replaced by a new almost complex structure such that new almost hermitian structure is of strict type $W_3\op W_4$, and belongs to the family of case (ii).
\end{enumerate}
The structure $(M^6,g,J)$ belongs to precisely one of these families.
\end{thm}

\begin{proof}
Most statements follow from Theorem \ref{thm_strict_types_and_isotropy_groups} and Proposition \ref{prop_TcscalarinCliff}. For case (iii), we note, that the isotropy group must be either $SU(2)$ or $T^1$, but that $SU(2)$ is not possible, Proposition \ref{prop_trivialholSU2impossible}.
\end{proof}

\noindent
The explicit classification of case (ii) is contained in Theorem \ref{thm_lie_groups_T2}.

%---------------------------------------------------------------------------- 
\subsubsection{Abelian holonomy}
%---------------------------------------------------------------------------- 
In the case of non-trivial, abelian holonomy, the dimension of the holonomy group must be either one or two. 

If the dimension is two, then $Hol_o(\nabla^c)=T^2_{-1}$, $T^2_1$ or $T^2=T^2_0$. Here $T^2_k$ denotes the maximal torus in $U(2)_k$. The structures are of strict type $W_1$, $W_3$ and $W_4$ respectively. In the first two cases the infinitesimal model is unique, and is realized on the twistor space $F(1,2)=SU(3)/T^2$. In the last case, the torsion form is unique, but the curvature depends on additional parameters. These structures, however, do not admit parallel spinors. 

If the dimension is one, there are infinitely many holonomy groups. If, however, we require the existence of a parallel spinor, then $Hol_o(\nabla^c)=T^1$. We discuss the classification in this case:

There are seven possible isotropy groups of the torsion form. It is known, that $Iso_o(T^c)=SU(3)$ or $U(2)_1$ is not possible. Likewise the case $Iso_o(T^c)=SO(3)$ can be excluded, Proposition \ref{lem_SO3_Hol_T1_impossible}. The remaining cases may be reduced to the case of structures related to $T^2$. This is obvious, if $Iso_o(T^c)=T^2$ or $U(2)_0$. If $Iso_o(T^c)=SU(2)$, the strict type is either $W_1\op W_3$ or $W_1\op W_3\op W_4$, and the almost complex structure may be replaced by a new structure such that the strict type is $W_4$ and $Iso(T^c)=U(2)_0$. If $Iso_o(T^c)=T^1$, the strict type is either $W_1\op W_3$ or $W_1\op W_3\op W_4$, and the almost complex structure may be replaced by a new structure such that the strict type is $W_3\op W_4$ and $T^2\subset Iso(T^c)$. The classification is now contained in Theorem \ref{thm_classification_holonomy_T1_T2}.

%---------------------------------------------------------------------------- 
\subsubsection{Non-abelian holonomy}
%---------------------------------------------------------------------------- 
In the case of non-abelian holonomy, the dimension of the holonomy group must be at least three. The six possible holonomy groups are $SU(3)$, $U(2)_{-1}$, $U(2)_1$, $U(2)_0$, $SU(2)$, and $SO(3)$.

If the strict type is $W_1$, the possible holonomy groups are $SU(3)$, $U(2)_{-1}$ and $SO(3)$. The infinitesimal model is unique in all cases, and is realized on $S^6$, $\CP^3$ and $S^3\x S^3$ respectively. Obviously, $Hol_o(\nabla^c)=SU(3)$ or $U(2)_{-1}$ implies strict type $W_1$. If $Hol_o(\nabla^c)=SO(3)$, strict type $W_3$ and strict type $W_1\op W_3$ is possible. In the case of strict type $W_3$, the infinitesimal model is unique, and is realized on $SL(2,\C)$. For the case of strict type $W_1\op W_3$, see the section on structures related to $SO(3)$.
 
If $Hol_o(\nabla^c)=U(2)_1$, the structure is of strict type $W_3$. The infinitesimal model is unique, and the structure may be realized on the twistor space $\CP^3$. It is related to the nearly K\"ahler on that space. 

If $Hol_o(\nabla^c)=U(2)_0$, the structure is of strict type $W_4$. The torsion form is unique, but there is one additional parameter for the curvature. These structures, however, do not admit parallel spinors. 

If $Hol_o(\nabla^c)=SU(2)_0$, the structure is of strict type $W_1\op W_3$, $W_1\op W_3\op W_4$ or $W_4$. The almost complex structure may be replaced by a new structure such that the strict type is $W_4$ and $Iso(T^c)=U(2)_0$.  The infinitesimal model in that case is unique, and is realized on $S^5\x S^1$ (cf. Example \ref{exa_generalized_Hopf_S5xS1}).

%\newpage
\bibliographystyle{amsalpha}
%

%
%%%%%%%%%%%%%%%%%%%%%%%%%%%%%%%%%%%%%%%%%%%%%%%%%%%%%%%%%%%%%%%%%%%%%%%%%%%%%
\end{document}